\documentclass[a4paper]{article}
      \usepackage[english]{babel}
  
      \usepackage[latin1]{inputenc}
      \usepackage{graphicx}
      \usepackage{multicol}
      \usepackage{graphics,epsfig}
      \usepackage{amsmath,amsfonts,amssymb}

      \usepackage{fancyhdr}

      \newtheorem{thm}{Theorem}[section]
      \newtheorem{propo}{Proposition}[section]
      \newtheorem{prop}{Properties}[section]
      \newtheorem{Def}{Definition}[section]
      \newtheorem{rmq}{Remark}[section]
      \newtheorem{lem}{Lemma}[section]
      \newtheorem{nota}{Notation}[section]
      
\title{\bf \large Wellposedness and regularity for a degenerate parabolic equation arising in a model of chemotaxis with nonlinear sensitivity}

 \author{\normalsize \bf Alexandre MONTARU\\
 \footnotesize Universit\'e Paris 13, Sorbonne Paris Cit\'e, \\
\footnotesize LAGA, CNRS, UMR 7539,\\ 
\footnotesize F-93430, Villetaneuse, France. \\
\small \textit{montaru@math.univ-paris13.fr} }

\date{}
\begin{document}
     
      \maketitle

\begin{abstract}
We study a one-dimensional parabolic PDE with degenerate diffusion and non-Lipschitz 
nonlinearity involving the derivative. This evolution equation arises when searching
 radially symmetric solutions of
a chemotaxis model of Patlak-Keller-Segel type. We prove its local in time wellposedness 
in 
some appropriate space, a blow-up alternative, regularity results and give an idea of the shape of solutions. A transformed and
an approximate problem naturally appear  in the way of the proof and are also  crucial in  \cite{Montaru} in order to study the global behaviour 
of solutions of the equation for a critical parameter, more precisely  to show the existence of a critical mass.  
\end{abstract}

\section*{Introduction}

In this paper, we are mainly interested in studying the local in time  wellposedness of 
the following problem $(PDE_m)$ :
\begin{eqnarray}
\label{equ_u_1}
u_t=x^{2-\frac{2}{N}}\,u_{xx}+u\,{u_x}^q & \qquad t>0 & 0<x\leq 1\\
\label{equ_u_2}
u(t,0)=0 &\qquad t\geq 0&\\
\label{equ_u_3}
u(t,1)=m &\qquad t\geq 0&\\
\label{equ_u_4}
u_x(t,x)\geq 0 &\qquad t> 0& 0\leq x\leq 1,
\end{eqnarray}
where $N$ is an integer greater or equal to $2$, $m\geq 0$ and $0<q<1$.\\

This problem follows from a chemotaxis model being aimed at describing a collection of 
cells diffusing and emitting a chemical which attracts themselves. These cells are  
assumed to lie in a physical domain corresponding to the open unit ball $D\subset \mathbb{R}^N$ ($N=2$ or $N=3$ being 
the most relevant cases) and if we suppose moreover that cells  diffuse 
much more slowly 
than the chemoattractant, we get the following  
parabolic-elliptic Patlak-Keller-Segel system $(PKS_q)$  :
\begin{eqnarray}
\rho_t= \Delta \rho-\nabla [\rho^q \, \nabla c]&\qquad t>0& \mbox{on } D\\
- \Delta c= \rho&\qquad t>0& \mbox{on } D
\end{eqnarray}
with the following boundary conditions :  
\begin{equation}
\label{BC-rho}
\frac{\partial \rho}{\partial \nu}-\rho^q \, \frac{\partial c}{\partial \nu}=0 \mbox{ 
\quad  on } \partial D
\end{equation}
\begin{equation}
\label{BC-c}
c=0 \mbox{\quad  on }\partial D 
\end{equation}

where $\rho$ is the cell density and $c$ the chemoattractant concentration. 
Note that on the boundary $\partial D$ are imposed a natural no flux condition  for 
$\rho$  
and Dirichlet conditions for $c$.\\
Problem $(PDE_m)$ follows from $(PKS_q)$ when considering radially symmetric solutions 
and after having made some transformations and a renormalization. \\
What is essential to 
know 
is that :
\begin{itemize}
\item $m$ is proportional to the cells mass $\int_B \rho $.
\item The derivative of $u$ is the quantity with physical interest since $u_x$ is 
proportional to the cells density $\rho$, up to a rescaling in time and a change of 
variable. More precisely, denoting $\rho(t,y)=\tilde{\rho}(t,|y|)$ for $t\geq 0$ and $y\in\overline{D}$,
$$ \tilde{\rho}(t,x)=N^{\frac{2}{q}}\,u_x(N^2\,t,x^N) \mbox{\qquad for all }x\in [0,1].$$
\item The power 
$q=\frac{2}{N}$ is critical.
\end{itemize}
Much more detail about problem $(PKS_q)$ and its link with $(PDE_m)$ are given in the 
introduction of \cite{Montaru}. See \cite{Patlak,KS,Herrero-Sastre} for references  concerning the biological background and \cite{Horstmann1,Horstmann2,Herrero,HP,Horstmann-Winkler,Kavallaris-Souplet,BKLN, BKLN2,Herrero-Velazquez,Dolbeault,Perthame,BCL,Blanchet2,Blanchet3,CALCOR} for related mathematical results.\\

The critical case $N=2,\,q=1$ is already well-known for its critical mass $8\pi$. See \cite{BKLN, 
Herrero-Velazquez}.
Our aim is to provide a rigorous framework  in view of the study that we have carried out
in \cite{Montaru} on the global behaviour of
solutions of problem $(PDE_m)$ in the case $N\geq 3$ and $q=\frac{2}{N}\in(0,1)$. In 
particular, we will prove the local in time existence  and uniqueness 
of a maximal  classical  solution $u$ for problem $(PDE_m)$ with initial 
condition $u_0\in Y_m$ where $Y_m$ is a space of functions which will be made explicit 
in the next section. Moreover, we have a blow-up alternative, regularity 
results and a description giving an idea of the shape of solutions.\\
Let us point out that solutions of $(PDE_m)$ are uniformly bounded in view of the maximum principle and that possible finite singularities are thus of gradient blow-up type. However, we shall show
(see Theorem \ref{thm_existence_u})iii)) that the solution can be continued as long as the slopes with respect to the origin are controlled, which is a crucial fact for the analysis in \cite{Montaru}.\\

In the way to prove these results, we will 
need some related problems, in particular a transformed problem $(tPDE_m)$ and an 
approximated problem $(PDE_m^\epsilon)$ for $\epsilon>0$.
We also would like to point out the role played by both  problems when proving in  
\cite{Montaru} that problem $(PDE_m)$ exhibits a critical mass 
phenomenon. More precisely, we 
showed there the existence of $M>0$ such that :
\begin{itemize}
\item If $m\leq M$, then $u$ is global and 
$$u(t)\underset{t\rightarrow \infty}{\longrightarrow} U \mbox{ \qquad in } 
C^1([0,1])$$
where $U$ is a steady state of $(PDE_m)$.
\item If $m>M$ then $u(t)$ blows up in finite time $T_{max}<\infty$. \\
Moreover,
$$\underset{t\rightarrow T_{max}}{\lim} \mathcal{N}[u(t)]=+\infty$$
 where 
$\mathcal{N}[f]=\underset{x\in(0,1]}{\sup}\frac{f(x)}{x}$ for any real function $f$ 
defined on $(0,1]$.

\end{itemize}

We precisely described the set of steady states and in particular proved that there 
exists 
only one stationary solution for $m<M$, none for $m>M$ but a whole continuum for $m=M$ 
(in which case $u_x$ has support strictly inside $[0,1))$. 
The critical case $m=M$ could then be much more intricate since the solution could for 
instance 
oscillate between various stationary solutions. In order to treat the case $m\leq M$, we 
used some dynamical systems methods and proved (with help of $(tPDE_m)$) that all trajectories 
are relatively compact and (with help of $(PDE^\epsilon_m)$)
the existence of a strict Lyapunov functional $\mathcal{F}=\underset{\epsilon \rightarrow 
0}{\lim}\,\mathcal{F}_\epsilon$ where $\mathcal{F}_\epsilon$ is a strict Lyapunov 
functional 
for $(PDE_m^\epsilon)$. \\

Eventually, we would like to stress that problem $(PDE_m)$ is not standard since it 
presents two difficulties :
\begin{itemize}
\item The diffusion is degenerate since $x^{2-\frac{2}{N}}$ goes to $0$ as $x$ goes to 
$0$.
\item The nonlinearity, which involves a gradient term, is not Lipschitz since $q\in(0,1)$.
\end{itemize}

The outline of the rest of the paper is as follows :
\tableofcontents

\section{Notation and strategy}

We give the definition of $Y_m$,  a space of functions appropriate for our study.
\begin{Def} Let $m\geq0$.\\
$Y_m=\{u\in C([0;1]),\:u \mbox{ nondecreasing },\; u'(0) \mbox{ exists, } u(0)=0,\;u(1)=m\}$
\end{Def}

We are interested in the following evolution equation called $(PDE_m)$ with  
$$N\geq 2,\; 
q\in(0,1) \text{ and } m\geq 0.$$ 

\begin{Def} Let $T>0$.\\
We define problem $(PDE_m)$ with initial condition $u_0\in Y_m$ by :
\begin{eqnarray} 
\label{equ_u1}
 & u_t=x^{2-\frac{2}{N}}u_{xx}+u \,{u_x}^q &\mbox{ on } (0,T]\times (0,1] \\
\label{equ_u2}
& u(0)=u_0 &\\
\label{equ_u3}
& u(t)\in Y_m  &\mbox{ for } t\in [0,T]
\end{eqnarray}
A classical solution of problem $(PDE_m)$ with initial condition $u_0\in 
Y_m$ on $[0,T]$   is a function
$$u\in C([0,T]\times [0,1])\bigcap
C^1((0,T]\times [0,1])\bigcap C^{1,2}((0,T]\times (0,1])$$
such that  
 $(\ref{equ_u1})(\ref{equ_u2})(\ref{equ_u3})$ are satisfied.\\
A classical solution of problem $(PDE_m)$ on $[0,T)$ is defined similarly.
\end{Def}

We would like to briefly describe the strategy used to obtain a maximal classical solution of 
problem $(PDE_m)$, as well as approximate solutions of it that turned out to be very helpful in
\cite{Montaru}. At the 
same time, we introduce the notation used throughout this paper.\\

\underline{First step :} we  introduce the change of unknown, denoted 
$\theta_0$, in order to get rid of the degenerate diffusion. It turns out
 (see formulae (\ref{u_i})(\ref{derivee_en_t_u_i})(\ref{derivee_u_i})(\ref{derivee_seconde_u_i})) 
 that the transformed equation becomes nondegenerate and involves
the radial heat operator, but in N+2 space dimensions.

\begin{Def} Let $B$ denote the open unit ball in 
$\mathbb{R}^{N+2} $. \\We define the transformation  \\
 $\begin{array}{ll}
\theta_0 : &Y_m\longrightarrow Z_m \\ 
&u\longrightarrow w \mbox{ where }w(y)= \frac{u(|y|^N)}{|y|^N}\mbox{ for all } y\in 
\overline{B}\backslash \{0\} 
\end{array}$\\

where $$Z_m=\{w\in C(\overline{B}), \; w| _ { \partial B } = m\}.$$

\end{Def}

\begin{rmq}
To avoid any confusion, we would like to stress that the physical domain  $D$ (where the cells live) lies in $\mathbb{R}^N$  but that the ball $B$ (where the transformed problem is posed) lies in $\mathbb{R}^{N+2}$.
\end{rmq}

Setting $w_0=\theta_0(u_0)\in Z_m$ and $w(t,y)=\frac{u(N^2t,|y|^N)}{|y|^N}$ for all $y\in 
\overline{B}$, we obtain a transformed problem 
called $(tPDE_m)$ with simple Laplacian diffusion which will allow us to use the heat semigroup.

\begin{Def}Let $m\geq 0$, and $T>0$. \\
Let $w_0=\theta_0(u_0)$ where $u_0\in Y_m$.\\
We define problem $(tPDE_m)$ with initial condition $w_0$ by :
\begin{eqnarray}
\label{equ_w1}
 & w_t=\Delta w + N^2 w \left(w+\frac{y.\nabla 
w}{N}\right)^q  &\mbox{ on }(0,T]\times \overline{B}\\
\label{equ_w2}
&w(0)=w_0&\\
\label{equ_w3}
&w+\frac{y.\nabla w}{N}\geq 0 & \mbox{ on }(0,T]\times \overline{B}\\
\label{equ_w4}
&w=m&\mbox{ on } [0,T]\times \partial B
\end{eqnarray}

 A classical solution of problem $(tPDE_m)$ with initial condition $w_0$ on $[0,T]$    is a function 
$$ w\in C([0,T]\times \overline{B}) \bigcap C^{1,2} ((0,T]\times \overline{B})$$
such that 
 $(\ref{equ_w1})(\ref{equ_w2})(\ref{equ_w3})(\ref{equ_w4})$ are satisfied.\\
We define analogously a classical solution on $[0,T)$.
\end{Def}

\underline{Second step :} since equation $(tPDE_m)$ still has a non Lipschitz nonlinearity, we want to 
define an approximate problem $(tPDE_m^\epsilon)$ for $\epsilon>0$ to get rid of it. \\
This is why we introduce the following function :

\begin{Def}
Let $\epsilon>0$. We set :
$$f_\epsilon(x)=(x+\epsilon)^q-\epsilon^q \mbox{ \quad if }x\geq 0 $$
and $f_\epsilon$ can be extended to $\mathbb{R}$ so that it satisfies both following conditions : 
$$f_\epsilon\in C^3(\mathbb{R})$$ 
$$-|x|^q\leq f_\epsilon(x)<0 \mbox{ for all } x\in (-\infty,0)$$
\end{Def}
Observe in particular that $|f_\epsilon(x)|\leq |x|^q$ for all $x\in\mathbb{R}$.

\begin{rmq}
Note that the conditions on $f_\epsilon$ on $(-\infty,0)$ are purely technical. Indeed,
 the choice of the 
extension does not play any role since we will
 prove that actually $(u^\epsilon)_x>0$ on $[0,1]$, where  $u^\epsilon$ is the maximal classical solution of problem $(PDE_m^\epsilon)$ with initial condition $u_0\in Y_m$ as defined below.
\end{rmq}

\begin{Def}Let $\epsilon>0$, $m\geq0$ and $T>0$. \\
Let $w_0=\theta_0(u_0)$ where $u_0\in Y_m$. \\
We define problem $(tPDE_m^\epsilon)$ with initial condition $w_0$ by :
\begin{eqnarray}
\label{equ_w_eps1}
 & w_t=\Delta w + N^2 w \;f_\epsilon\left(w+\frac{y.\nabla 
w}{N}\right)  &\mbox{ on }(0,T]\times \overline{B}\\
\label{equ_w_eps2}
&w(0)=w_0\\
\label{equ_w_eps3}
&w=m&\mbox{ on } [0,T]\times \partial B
\end{eqnarray}
 A classical solution for problem $(tPDE_m^\epsilon)$ with initial condition $w_0$ on $[0,T]$ is a function 
$$ w^\epsilon\in C([0,T]\times \overline{B}) \bigcap C^{1,2} ((0,T]\times \overline{B})$$
such that
 $(\ref{equ_w_eps1})(\ref{equ_w_eps2})(\ref{equ_w_eps3})$ are satisfied.\\
We define similarly a classical solution on $[0,T)$.
\end{Def}

The setting of problem $(tPDE_m^\epsilon)$ is standard and allows to find a 
unique classical maximal solution $w^\epsilon$  on $[0,T_\epsilon^*)$ with initial 
condition $w_0=\theta(u_0)$ for any $u_0\in Y_m$.
Then, a compactness property and the monotonicity  of the family 
$(w^\epsilon)_{\epsilon>0}$ 
allows to get a local solution of $(tPDE_m)$ by letting $\epsilon$ go to 
$0$. Eventually, since a 
comparison principle is available, we obtain a unique maximal classical 
solution $w$ for problem 
$(tPDE_m)$.
Since $w_0$ is radial, so is $w(t)$ which can then be written 
$w(t,y)=\tilde{w}(t,|y|)$ for all $y\in \overline{B}$.
Eventually, setting $$u(t,x)=\tilde{w}(\frac{t}{N^2},x^{\frac{1}{N}})\mbox{\qquad for all }
x\in[0,1],$$ we  get a 
classical solution for problem $(PDE_m)$ that will be proved to be actually maximal.\\

As explained before, we  will also need solutions of $(PDE_m^\epsilon)$, an approximate 
version  of problem $(PDE_m)$.

\begin{Def} Let $\epsilon>0$, $m\geq 0$ and $T>0$.\\
 We define problem $(PDE_m^\epsilon)$ with initial condition $u_0\in Y_m$ by :
\begin{eqnarray} 
\label{equ_u_eps1}
 & u_t=x^{2-\frac{2}{N}}u_{xx}+u f_\epsilon(u_x) &\mbox{ on } (0,T]\times (0,1] \\
\label{equ_u_eps2}
& u(0)=u_0 & \\
\label{equ_u_eps3}
&u(t,0)=0 &\mbox{ for all } t\in [0,T]\\
\label{equ_u_eps4}
&u(t,1)=m &\mbox{ for all } t\in [0,T]
\end{eqnarray} 

 A classical solution of problem $(PDE_m^\epsilon)$ with initial condition
 $u_0\in Y_m$ on $[0,T]$  is a function 
$$u^\epsilon\in C([0,T]\times [0,1])\bigcap
C^1((0,T]\times [0,1])\bigcap C^{1,2}((0,T]\times (0,1])$$
such that
 $(\ref{equ_u_eps1})(\ref{equ_u_eps2})(\ref{equ_u_eps3})(\ref{equ_u_eps4})$ are satisfied.\\
A classical solution of problem $(PDE_m^\epsilon)$ on $[0,T)$ is defined similarly.
\end{Def}

We will see that each of the four problems we have described admits a unique maximal classical solution and we would like to fix now the notation we will use throughout this paper for these solutions.
\begin{nota}\quad
\begin{itemize}
\item Let $u_0\in Y_m$.\\
We denote $u$ [resp. $u^\epsilon$] the maximal classical solution of problem $(PDE_m)$ [resp. $(PDE_m^\epsilon)$] with initial condition $u_0$.
\item Let $w_0=\theta_0(u_0)$ where $u_0\in Y_m$.\\
We denote $w$ [resp. $w^\epsilon$] the maximal classical solution of problem $(tPDE_m)$ [resp. $(tPDE_m^\epsilon)$] with initial condition $w_0$.
\end{itemize}
\end{nota}

\section{Main results : local wellposedness, regularity and blow-up alternative for problem $(PDE_m)$}

\begin{Def} For any real function $f$ defined on $(0,1]$, we set
$$\mathcal{N}[f]=\underset{x\in(0,1]}{\sup}\frac{f(x)}{x}.$$
\end{Def}

\begin{thm}Let $N\geq 2$, $q\in(0,1)$ and $m\geq 0$. \\
Let  $K>0$ and $u_0\in Y_m $ with $\mathcal{N}[u_0]\leq K$.
\label{thm_existence_u}
\begin{itemize}
 \item [i)]
There exists $T_{max}=T_{max}(u_0)>0$ and a unique maximal classical solution $u$ of problem $(PDE_m)$ with initial condition $u_0$. \\
Moreover, $u$ satisfies the following condition :
$$ \underset{t\in(0,T]}{\sup} \sqrt{t}\; \|u(t)\|_{C^1([0,1])}<\infty \mbox{ for any } T\in(0,T_{max})$$
\item [ii)] There exists $\tau=\tau(K)>0$ such that $T_{max}\geq \tau$.
\item[iii)] Blow up alternative : 
$T_{max}=+\infty$ \quad or \quad  $\underset{t\rightarrow T_{max}}
{\lim}\mathcal{N}[u(t)]=+\infty $

\item[iv)] $u_x(t,0)>0$ for all $t\in (0,T_{max})$.
\item [v)]
If $0<t_0<T<T_{max}$ and $x_0\in(0,1)$, then for any $\gamma \in(0,q)$,  
$$u\in C^{1+\frac{\gamma}{2},2+\gamma} ( [t_0,T]\times [x_0,1])$$
\item [vi)] For all $t\in(0,T_{max})$, $u(t)\in Y_m^{1,\frac{2}{N}}$ where for any 
$\gamma>0$,
$$Y_m^{1,\gamma}=\{u\in Y_m\cap C^1([0,1]), \;\underset{x\in (0,1]}{\sup}\; \frac{|u'(x)-u'(0)|}{x^\gamma}<\infty\}.$$
\end{itemize}
\end{thm}

Remember that the radially symmetric cells density $\rho$ is related to the derivative of $u$ by :
$$ \tilde{\rho}(t,x)=N^{\frac{2}{q}}\,u_x(N^2\,t,x^N) \mbox{\qquad for all }x\in [0,1].$$
We can have an idea of the shape of $u_x$, especially near the origin since we can show :

\begin{propo}
\label{propo_forme} Let $u_0\in Y_m$.
\begin{itemize} 
\item[i)] For all $(t,x)\in (0,T_{max}(u_0))\times [0,1]$,  $$u_x(t,x)=h(t,x^{\frac{1}{N}})$$
 with $h\in C^{1,1}((0,T_{max}(u_0))\times [0,1])$.
\item[ii)] Let $[t_0,T]\subset (0,T_{max}(u_0))$. \\
There exists $\delta>0$ such 
that for all $(t,x)\in [t_0,T]\times [0,\delta]$,
$$u_x(t,x)=h(t,x^{\frac{1}{N}})$$
 with $h\in C^{1,\infty}([t_0,T]\times [0,1])$ such that for any $t\in[t_0,T]$, \\ $h(t,\cdot)$ has odd derivatives vanishing at $x=0$.
\item[iii)] Let $t\in(0,T_{max}(u_0))$.\\
$u_x(t,x)$ admits an expansion of any order in powers of $x^\frac{2}{N}$ at $x=0$.\\
For instance, $u_x(t,x)=a(t)+b(t)x^\frac{2}{N}+o(x^\frac{2}{N})$.
\end{itemize}
\end{propo}

\section{Additional results}

\subsection{Problem $(tPDE_m)$}

\begin{thm} 
\label{thm_existence_w} 
Let $u_0\in Y_m$ and $w_0=\theta_0(u_0)$.
\begin{itemize}
 \item [i)] There exists $T^*=T^*(w_0)>0$ and a unique maximal 
classical solution $w$ of problem $(tPDE)$ with initial condition $w_0$.\\
Moreover, $w$ satisfies the following condition :
$$ \underset{t\in(0,T]}{\sup} \sqrt{t}\; \|w(t)\|_{C^1(\overline{B})}<\infty \mbox{ for any } T\in(0,T^*).$$
\item [ii)] Blow-up alternative : $T^*=+\infty \quad$ or $\quad\underset{t\rightarrow T^* }{\lim}
\|w(t)\|_{\infty, \overline{B}}=+\infty$.
 \item[iii)] $w>0$ on $(0,T^*)\times \overline{B} $.
\item[iv)] $w\in C^{1+\frac{\gamma}{2},2+\gamma}([t_0,T]\times \overline{B})$ for all $\gamma\in (0,q)$ and all $[t_0,T]\subset (0,T^*)$.
\end{itemize}
\end{thm}

\textbf{Connection with problem $(PDE_m)$ :} 
$$T_{max}(u_0)=N^2 T^*(w_0)$$ 
and 
for all $(t,x)\in [0,T_{max})\times [0,1]$, 
\begin{equation}
\label{def_u}
u(t,x)=x\; \tilde{w}(\frac{t}{N^2},x^{\frac{1}{N}})
\end{equation}
 where for any radially symmetric function $f$ on $B$, we will denote $f(y)=\tilde{f}(|y|)$ for all $y\in B$.

\subsection{Problem $(PDE_m^\epsilon)$}

\begin{thm}Let $m\geq 0$, $\epsilon>0$ and $K>0$. 
\\Let $u_0\in Y_m $ with $\mathcal{N}[u_0]\leq K$.
\label{thm_existence_u_epsilon}
\begin{itemize}
 \item [i)]
There exists $T_{max}^\epsilon=T_{max}^\epsilon(u_0)>0$ and a unique maximal 
classical solution $u^\epsilon$ on $[0,T_{max}^\epsilon)$ of  problem
$(PDE_m^\epsilon)$ with initial condition $u_0$.\\
Moreover, $u^\epsilon$ satisfies the following condition :
\begin{equation}
\label{maj_u_eps}
 \underset{t\in(0,T]}{\sup} \sqrt{t}\; \|u^\epsilon(t)\|_{C^1([0,1])}<\infty \mbox{ for all } T\in(0,T_{max}^\epsilon).
\end{equation}
\item [ii)] There exists $\tau=\tau(K)>0$ such that for all $\epsilon>0$, 
$T_{max}^\epsilon\geq \tau$. \\
Moreover, there exists $C=C(K)>0$  independent of $\epsilon$ such that
\begin{equation}
\label{majoration_u_epsilon_tau'}
 \underset{t\in [0,\tau]}{\sup} \mathcal{N}[u^\epsilon(t)] 
\leq C.
\end{equation} 
\item[iii)] Blow up alternative : 
$T^\epsilon_{max}=\infty$ \quad or \quad  $\underset{t\rightarrow T_{max}^\epsilon}
{\lim}\mathcal{N}[u^\epsilon(t)]=\infty $.
\item[iv)] $(u^\epsilon)_x>0$ on $ (0,T_{max}^\epsilon)\times [0,1]$.
\item [v)]
If $0<t_0<T<T_{max}^\epsilon$ and $x_0\in(0,1)$, then for any $\gamma \in(0,1)$,  
$$u^\epsilon\in C^{1+\frac{\gamma}{2},2+\gamma} ( [t_0,T]\times [x_0,1]).$$
\item[vi)] If $u_0\in Y_m^{1,\gamma}$ with $\gamma>\frac{1}{N}$ then $u^\epsilon\in 
C([0,T_{max}^\epsilon),C^1([0,1]))$.
\end{itemize}
\end{thm}

\textbf{Connection with problem $(PDE_m)$ :} \\
Fixing an initial condition $u_0\in Y_m$, 
the next lemma shows the convergence of maximal classical solutions $u^\epsilon$ of 
$(PDE_m^\epsilon)$ to the maximal classical solution of $(PDE_m)$ in various spaces. \\
These results turned out to be essential in \cite{Montaru} since, starting from a strict
 Lyapounov functional $\mathcal{F}_\epsilon$ for 
$(PDE_m^\epsilon)$ in the subcritical case ($m$ less or equal to the critical mass $M$), 
we obtained a strict Lyapounov functional $\mathcal{F}$ for $(PDE_m)$ by setting 
$\mathcal{F}=\underset{\epsilon \rightarrow 0}{\lim}\,\mathcal{F}_\epsilon$.
We point out that it does not seem possible to construct a Lyapunov functional for $(PDE_m)$ by a direct approach (cf. p.7 in \cite{Montaru}).

\begin{lem}
\label{lem_convergence_u_epsilon}
Let $u_0\in Y_m$.
\begin{itemize}
\item[i)] $T_{max}(u_0)\leq T_{max}^\epsilon(u_0)$ for any $\epsilon>0$.
\item[ii)]  Let $[t_0,T]\subset (0,T_{max}(u_0))$.
\begin{itemize}
\item [$\alpha)$]$u^\epsilon \underset{\epsilon \rightarrow 0}{\longrightarrow} u$ in $C^{1,2}([t_0,T]\times (0,1]).$\\
Moreover, there exists $K>0$ independent of $\epsilon$ such that\\
for all $(t,x)\in[t_0,T]\times (0,1]$, $|u^\epsilon_{xx}|\leq \frac{K}{x^{1-q}}$.
\item[$\beta )$] $(u^\epsilon)_x \underset{\epsilon \rightarrow 0}{\longrightarrow} u_x$ in $C([t_0,T]\times [0,1])$.
\item[$\gamma )$] $(u^\epsilon)_t \underset{\epsilon \rightarrow 0}{\longrightarrow} u_t$ in $C([t_0,T]\times [0,1])$.
\end{itemize}
\end{itemize} 
\end{lem}

\textbf{Connection with problem $(tPDE_m^\epsilon)$.}\\
 Let $u_0\in Y_m$ and $w_0=\theta_0(u_0)$. Then 
 $$T^\epsilon_{max}(u_0)=N^2 T_\epsilon^*(w_0).$$
Moreover, for all $(t,x)\in [0,T_{max}^\epsilon)\times [0,1]$, $$u^\epsilon(t,x)=x\; \tilde{w}^\epsilon(\frac{t}{N^2},x^{\frac{1}{N}}).$$

\section{Proofs}

\subsection{Comparison principles}

The four problems we have defined each admit a comparison principle which is in 
particular available  
for classical solutions. \\
Whence the uniqueness of the maximal classical solution in each case.

\begin{lem}
\label{PC_PDE}
 \textbf{Comparison principle for problem $(PDE_m)$} \\
Let $T>0$. Assume that :
\begin{itemize}
\item $u_1,u_2\in C([0,T]\times [0,1])\bigcap
C^1((0,T]\times [0,1])\bigcap C^{1,2}((0,T]\times (0,1))$. 
\item For all $t\in (0,T] $, $u_1(t)$ and $u_2(t)$ are nondecreasing. 
\item There exists $i_0\in\{1,2\}$ and some $\gamma<\frac{1}{q}$ such that : 
$$ \underset{t\in(0,T]}{\sup} t^\gamma\; \|u_{i_0}(t)\|_{C^1([0,1])}<\infty .$$
\end{itemize}
Suppose moreover that :
\begin{eqnarray} 
 & (u_1)_t\leq x^{2-\frac{2}{N}}(u_1)_{xx}+u_1 (u_1)_x^q  & \mbox{ for all }(t,x)\in (0,T]\times (0,1).\\
 & (u_2)_t\geq x^{2-\frac{2}{N}}(u_2)_{xx}+u_2 (u_2)_x^q  & \mbox{ for all }(t,x)\in (0,T]\times (0,1).\\
& u_1(0,x)\leq u_2(0,x) &\mbox{ for all }x\in [0,1].\\
& u_1(t,0)\leq u_2(t,0) &\mbox{ for }t\geq 0. \\
& u_1(t,1)\leq u_2(t,1) &\mbox{ for }t\geq 0.
\end{eqnarray}
Then $u_1\leq u_2$ on $[0,T]\times [0,1]$.\\
\end{lem}

\textit{Proof :} 
Let us set $z=(u_1-u_2)e^{-\int_0^t (\|u_{i_0}(s)\|^q_{C^1}+1)ds}$. The hypotheses made show that
$z\in C([0;T]\times [0;1])\bigcap C^1((0,T]\times [0,1])\bigcap C^{1,2}((0,T]\times (0,1))$. \\
Assume now by contradiction that $\underset{[0;T]\times [0;1]}{\max} z>0$. \\
By assumption, $z\leq 0$ on the parabolic boundary of $[0,T]\times [0,1]$. \\
Hence,
$\underset{[0;T]\times [0;1]}{\max} z$ is reached at a point $(t_0,x_0)\in (0;T]\times (0;1)$. \\Then $z_x(t_0,x_0)=0$ so  
$(u_1)_x(t_0,x_0)=(u_2)_x(t_0,x_0)$. \\
Moreover, $z_{xx}(t_0,x_0)\leq 0$ and $z_t(t_0,x_0)\geq 0$.\\
 We have $z_t(t_0,x_0)\leq x^{2-\frac{2}{N}}z_{xx}(t_0,x_0)+ \left[ (u_{i_0})_x(t_0,x_0)^q-\|u_{i_0}(t_0)\|^q_{C^1}-1\right] z(t_0,x_0)$.
 The LHS of the inequality is nonnegative and the RHS is negative, whence the 
 contradiction.
 
\begin{rmq} \textbf{Comparison principle for problem $(PDE^\epsilon_m)$} \\
 Under the same assumptions (except the monotonicity of $u_1(t)$ and $u_2(t)$), an analogous comparison principle is available for problem $(PDE_m^\epsilon)$ for any $\epsilon>0$.

\end{rmq}

\begin{lem} \textbf{Comparison principle for problem $(tPDE_m)$} \\
Let $T>0$. Assume that :
\begin{itemize}
\item $w_1,w_2\in C([0,T]\times [0,1])\bigcap C^{1,2}((0,T]\times \overline{B})$. 
\item For $i=1,2$, for all $(t,y)\in (0,T]\times \overline{B} $, $w_i(t,y)=\tilde{w}_i(t,|y|).$
\item For $i=1,2$, for all $(t,y)\in (0,T]\times \overline{B} $, $w_i(t,y)+\frac{y.\nabla w_i(t,y)}{N}\geq 0$.
\item There exists $i_0\in\{1,2\}$ and some $\gamma<\frac{1}{q}$ such that :
 $$ \underset{t\in(0,T]}{\sup} t^\gamma\; \|\tilde{w}_{i_0}(t)\|_{C^1([0,1])}<\infty .$$
\end{itemize}
Suppose moreover that :
\begin{eqnarray} 
 & (w_1)_t\leq  \Delta w_1 + N^2 w_1 (w_1+\frac{y.\nabla w_1}{N})^q  &\mbox{ on }(0,T]\times \overline{B}.\\
 & (w_2)_t\geq  \Delta w_2 + N^2 w_2 (w_2+\frac{y.\nabla w_2}{N})^q  &\mbox{ on }(0,T]\times \overline{B}.\\
& w_1(0,y)\leq w_2(0,y) &\mbox{ for all } y\in \overline{B}.\\
& w_1(t,y)\leq w_2(t,y) &\mbox{ for all } (t,y)\in [0,T]\times \partial{B}.\\
\end{eqnarray}
Then $w_1\leq w_2$ on $[0,T]\times \overline{B}$.

\end{lem}
\textit{Proof :} For $i=1,2$, let us set 
\begin{equation}
\label{u_i}
u_i(t,x)=x\; \tilde{w}_i(\frac{t}{N^2},x^{\frac{1}{N}}).
\end{equation}
Calculations show that, for $0<t\leq T$ and $0<x\leq 1$ :
\begin{eqnarray}
 \label{derivee_en_t_u_i}
(u_i)_t(t,x)&=&\frac{x}{N^2} \; (\tilde{w}_i)_{t}
(\frac{t}{N^2},x^{\frac{1}{N}}).\\
 \label{derivee_u_i}
(u_i)_x(t,x)&=&\left[\tilde{w}_i+\frac{r(\tilde{w}_i)_r}{N}\right] 
(\frac{t}{N^2},x^{\frac{1}{N}})\nonumber \\
&=&\left[w_i+\frac{y.\nabla w_i)}{N}\right] 
(\frac{t}{N^2},x^{\frac{1}{N}}).\\
 \label{derivee_seconde_u_i}
x^{2-\frac{2}{N}}(u_i)_{xx}(t,x)&=&\frac{x}{N^2}\left[(\tilde{w}_i)_{rr}+\frac{N+1}{r}(\tilde{w}_i)_{r}\right] 
(\frac{t}{N^2},x^{\frac{1}{N}})\nonumber\\
&=&\frac{x}{N^2}\; \Delta w_i (\frac{t}{N^2},x^{\frac{1}{N}}).
\end{eqnarray}
 It is easy to check that $$u_i\in C([0,N^2T]\times 
 [0,1])\bigcap C^1((0,N^2T]\times [0,1])\bigcap
 C^{1,2}((0,N^2T]\times (0,1]).$$
Special attention has to be paid to the fact that $u_i$ is $C^1$ up to $x=0$ but this is 
clear because of (\ref{u_i}) and (\ref{derivee_u_i}).\\
Clearly, $u_1$ and $u_2$ satisfy all assumptions of Lemma \ref{PC_PDE}, so $u_1\leq u_2$ on $[0,T]
\times [0,1]$. Then $w_1\leq w_2$ on  $[0,T]\times \overline{B}\backslash \{0\}$. But by continuity of $w_1$ and $w_2$, we get $w_1\leq w_2$ on $[0,T]\times \overline{B}$.

\begin{rmq}
A similar comparison principle is available for problem $(tPDE_m^\epsilon)$ for any $\epsilon>0$ (except that we do not have to suppose $w_i(t,y)+\frac{y.\nabla w_i(t,y)}{N}\geq 0$ for $i=1,2$).
\end{rmq}

\subsection{Preliminaries to local existence results}

First, we would like to recall some notation and properties of the heat semigroup.
For reference, see for instance the book \cite{Lu} of A. Lunardi.

\begin{nota} \quad 
\begin{itemize}
\item $B$ denotes the open unit ball in $\mathbb{R}^{N+2}$.
\item $X_0=\{W \in C(\overline{B}), \; 
W | _{\partial B}=0 \}$.
\item $(S(t))_{t\geq 0}$ denotes the heat semigroup  on $X_0$. It is the restriction on 
$X_0$ of the Dirichlet heat semigroup on $L^2(B)$.
\item $(X_\theta)_{\theta \in [0;1]}$ denotes the scale of interpolation spaces for 
$(S(t))_{t\geq 0}$.
\end{itemize}

\end{nota}

\begin{prop}\quad 
\begin{itemize}
\item $X_{\frac{1}{2}}=\{ W \in C^1(\overline{B}), \;W | _{\partial B}=0 \}.  $
\item Let $\gamma_0 \in (0;\frac{1}{2}]$. For any $\gamma\in[0,2\gamma_0)$,  $$X_{\frac{1}{2}+\gamma_0} \subset C^{1,\gamma}
(\overline{B})$$  with continuous embedding.
\item There exists $C_D\geq 1$ such that for any $\theta \in [0;1]$,
 $W \in C(\overline{B})$ and $t>0$,
$$\|S(t)W \|_{X_\theta}\leq \frac{C_D}{t^{\theta}}\|W\|_{\infty}. $$
\end{itemize}

\end{prop}

For reference, we recall some notation and then introduce two spaces of functions more in order to state a useful lemma on $\theta_0$.
\begin{nota} Let $m\geq 0$ and $\gamma>0$. 
\begin{itemize}
\item For $W\in C^1(\overline{B})$, the $C^1$ norm of $W$ is $\| W\|_{C^1}=\|
W\|_{\infty,\overline{B}}+\|\nabla W\|_{\infty,\overline{B}}$.
\item $Y_m=\{u\in C([0;1]) \mbox{ nondecreasing},\; u'(0) \mbox{ exists, } u(0)=0,u(1)=m\}$.
\item $Z_m=\{w\in C(\overline{B}), \; w| _ { \partial B } = m\}$.
\item  $Y_m^{1,\gamma}=\{u\in Y_m\cap C^1([0,1]), \;\underset{x\in (0,1]}{\sup}\; \frac{|u'(x)-u'(0)|}{x^\gamma}<\infty\}$.
\item $Z_m^{1,\gamma}=\{w\in Z_m\cap C^1(\overline{B}),\; 
\underset{y \in \overline{B}\backslash \{0\}}{\sup}\; \frac{|\nabla w(y)|}{ |y|^{\gamma}}<\infty \}$ .
\item $\begin{array}{ll}
\theta_0 : &Y_m\longrightarrow Z_m \\ 
&u\longrightarrow w \mbox{ where }w(y)= \frac{u(|y|^N)}{|y|^N}\mbox{ for all } y\in 
\overline{B}\backslash \{0\}, \,w \mbox{ continuous on } \overline{B}. 
\end{array}$
\item Let $(a,b)\in(0,1)^2$. We denote $I(a,b)=\int_0^1 \frac{ds}{(1-s)^as^b}$.  \\ For all $ t\geq 0,\;
\int_0^t \frac{ds}{(t-s)^as^b}=t ^{1-a-b}I(a,b) $.
\end{itemize}

\end{nota}

\begin{lem} 
\label{propriete_theta}
 Let $m\geq 0$.
\begin{itemize}

\item[i)] $\theta_0$ sends $Y_m$ into $Z_m$.
\item [ii)]Let $\gamma>\frac{1}{N}$. 
$\theta_0$ sends $Y_m^{1,\gamma}$ into $Z_m^{1,N\gamma-1 }$.
\end{itemize}
\end{lem} 

\textit{Proof :} i) Let $u\in Y_m$ and $w=\theta_0(u)$. Clearly, $w$ can be extended in a continuous function on $\overline{B}$ by setting $w(0)=u'(0)$.\\
ii) Let $u\in Y_m^{1,\gamma}$ . 
It is clear that $w\in C^1(\overline{B}\backslash \{0\})$.\\
Let $y\in \overline{B}\backslash \{0\}$. 
$w(y)=\int_0^1 u'(t|y|^N)dt =w(0)+\int_0^1 [u'(t|y|^N)-u'(0)]dt$. \\
 Since $u\in Y_m^{1,\gamma}$,  there exists $K>0$ such that $|w(y)-w(0)|\leq K |y|^{N\gamma}$. Since $N\gamma>1$, $w$ is differentiable at $y=0$ and $\nabla w(0)=0$.\\
$\nabla w(y)=N\frac{y}{|y|^2}[u'(|y|^N)-w(y)]=
N\frac{y}{|y|^2}[u'(|y|^N)-u'(0)+w(0)-w(y)]$
So  $|\nabla w(y)| \leq 2NK |y|^{N\gamma-1}$. \\
Then $w\in C^1(\overline{B})$ and $\underset{y \in \overline{B}\backslash \{0\}}{\sup}\; \frac{|\nabla w(y)|}{ |y|^{N\gamma-1}}<\infty $, ie $w\in Z_m^{1,N\gamma-1 }$.\\

\begin{lem} \textbf{A density lemma.}\\
\label{lem_density}
Let $u_0\in Y_m$.
There exists a sequence $(u_n)\in Y_m^{1,1}$  such that $$\|u_n-u_0\|_{\infty,[0,1]} \underset{n\rightarrow \infty}{\longrightarrow} 0$$ and
$$\mathcal{N}(u_n)\leq \mathcal{N}(u_0).$$
\end{lem}

\textit{Proof :} Let $\epsilon >0$. Let $\mathcal{T}=\{(x,y) \in \mathbb{R}^2, \; 0\leq x \leq 1, \; 0\leq y\leq \mathcal{N}[u_0] \, x \}$. The graph 
of $u_0$ lies inside $\mathcal{T}$. Since $u_0$ is uniformly continuous on $[0,1]$, for $n_0$ large enough, we can construct a nondecreasing piecewise affine function on $[0,1]$ $v\in Y_m$  such that $\|u_0-v\|_{\infty,[0,1]}\leq \frac{\epsilon}{2}$ and for all $k=0\dots n_0$, $v(x_k)=u_0(x_k)$ where
$x_k=\frac{k}{n_0}$ and $v$ is affine between the successive points 
$P_k=(x_k,u(x_k))$ .  \\
Since $\mathcal{T}$ is convex and all points $P_k$ are in $\mathcal{T}$ 
then the graph of $v$ lies also inside $\mathcal{T}$.\\
We now just have to find a function $w\in C^2([0,1])\bigcap Y_m$ whose graph is in 
$\mathcal{T}$ and such that 
$\|w-v\|_{\infty,[0,1]}\leq \frac{\epsilon}{2}$.\\
In order to do that, we extend $v$ to a nondecreasing function $\overline{v}$ on 
$\mathbb{R}$ : we simply extend the first and last segments $[P_0,P_1]$ and 
$[P_{n_0-1},P_{n_0}]$ to a straight line, so that $\overline{v}$ is in 
particular affine on $(-\infty,\frac{1}{n_0}]$ and $[1-\frac{1}{n_0},+\infty)$. \\
Let $(\rho_\alpha)_{\alpha>0}$ a mollifiers family such that  
$\int_\mathbb{R}\rho_\alpha=1$, $supp(\rho_\alpha)\subset [-\alpha,
\alpha]$ and $\rho_\alpha$ is even.\\ 
Since $\overline{v}$ is Lipschitz continuous, there exists $\alpha_0>0$ such for all
$(x,y)\in\mathbb{R}^2$, if $|x-y|\leq \alpha_0$ then $|\overline{v}
(x)-\overline{v}(y)|\leq \frac{\epsilon}{2}$.  
Let $\alpha_1=\min(\alpha_0,\frac{1}{2n_0})$ and 
$w=\rho_{\alpha_1}*\, \overline{v}$. 
Remark that since $v$ is nondecreasing, so is $w$   and
$\|w-v\|_{\infty,[0,1]}\leq \frac{\epsilon}{2}$. \\
Note that, since $\int_\mathbb{R} 
 y\rho_{\alpha_1}(y)dy=0$, if for all
 $y\in [x-\alpha_1,x+\alpha_1]$, 
$\overline{v}(y)=ay+b$ (resp. $\overline{v}(y)\leq ay+b$)  then $w(x)=ax+b$ (resp. $w(x)\leq ax+b$).\\
Since the graph of $v$ lies inside $\mathcal{T}$ on $[0,1]$, this implies that
the graph of $w$ lies inside $\mathcal{T}$ on $[\dfrac{1}{2n_0},1-\dfrac{1}{2n_0}]$.\\
Moreover, since $\overline{v}$ is affine on $(-\infty,\frac{1}{n_0}]$ and on    
$[1-\frac{1}{n_0},+\infty)$, then $w$ is affine and coincides with $v$ on $[0,\frac{1}{2n_0}]$ and on 
$[1-\frac{1}{2n_0},1]$. \\
So $w\in Y_m$ and the graph of $w$ on $[0,1]$ lies inside $\mathcal{T}$.\\
Finally, $\|w-u_0\|_{\infty,[0,1]}\leq \epsilon$ and $w\in C^2([0,1])\bigcap Y_m\subset Y_m^{1,1}$.

\subsection{Solutions of problem ($tPDE_m^\epsilon$)}

\begin{thm} \textbf{Wellposedness of problem $(tPDE_m^\epsilon)$.}\\
Let $\epsilon>0$ and $K>0$.\\
Let $w_0\in Z_m $ with $\|w_0\|_{\infty,\overline{B}}\leq K$.
\label{thm_existence_w_epsilon}
\begin{itemize}
 \item [i)]
There exists $T^*_\epsilon=T^*_\epsilon(w_0)>0$ and a unique maximal classical solution 
$w^\epsilon$ of problem $(tPDE_m^\epsilon)$ on $[0,T_\epsilon^*)$ with initial condition 
$w_0$. \\
Moreover, $w^\epsilon$ satisfies the following condition :
\begin{equation}
\label{maj_w_eps}
 \underset{t\in(0,T]}{\sup} \sqrt{t}\; \|w^\epsilon(t)\|_{C^1(\overline{B})}<\infty 
 \mbox{ for all } T\in(0,T^*_\epsilon).
\end{equation}
\item [ii)]
We have the following blow-up alternative  :
$$ T^*_\epsilon=+\infty \quad \mbox{ or }\quad  \underset{t \rightarrow T^*_\epsilon}
{\lim} \|w^\epsilon(t) \|
_{\infty,\overline{B}}=+\infty.$$
\item [iii)] There exists $\tau=\tau(K)>0$ such that for all $\epsilon>0$,   
$T^*_\epsilon\geq \tau$.\\
Moreover, there exists $C=C(K)>0$  independent of $\epsilon$ such that
$$\underset{t \in [0,\tau]}{\sup}\|w^\epsilon(t)\|_{\infty, \overline{B}}\leq C . $$
\item[iv)]There exist $\tau'=\tau'(K)>0$ 
and $C'=C'(K)>0$  both independent of $\epsilon$ such that
$$\underset{t \in (0,\tau']}{\sup} \sqrt{t}\|w^\epsilon(t)\|_{C^1(\overline{B})}\leq 
C'  .$$
\item [v)]
If $0<t_0<T<T^*_\epsilon$, then for any $\gamma \in(0,1)$,  
$$w^\epsilon\in C^{1+\frac{\gamma}{2},2+\gamma} ( [t_0,T]\times \overline{B}).$$
\item [vi)] $w^\epsilon >0$ on $(0,T^*_{\epsilon'})\times \overline{B}$.  
\item [vii)] If $w_0\in Z_m\bigcap C^1(\overline{B})$, then $w^\epsilon \in 
C([0,T_\epsilon^*),C^1(\overline{B}))$.
\end{itemize}
\end{thm}

The proof of this theorem is based on a series of lemmas. We start with the following 
small time existence result for the auxiliary problem obtained by setting $W=w-m$ in 
$(tPDE_m^\epsilon)$.

\begin{lem}
 \label{lem_w_app}
Let $m\geq 0$, $\epsilon>0$ and $W_0 \in X_0 $.\\
There exists $\tau=\tau(\epsilon)>0$ and a unique mild solution 
$$W^\epsilon \in C([0;\tau], X_0) \bigcap
\{ W \in L^\infty_{loc}((0;\tau], X_{\frac{1}{2}}), \; \underset{ t\in (0;\tau] }{\sup} 
\sqrt{t} \;\| W(t) \| _ {C^1} < \infty  \}$$  of the following problem :
\begin{eqnarray} 
\label{pb_auxiliaire_1}
  W_t=\Delta W + N^2 (m+W) f_\epsilon(m+W+\frac{y.\nabla 
W}{N})  &\mbox{ on }&(0,\tau]\times \overline{B}\\
\label{pb_auxiliaire_2}
 W=0&\mbox{ on }&[0,\tau]\times \partial B\\
 \label{pb_auxiliaire_3}
 W(0)=W_0&&
\end{eqnarray}
More precisely, $W^\epsilon\in C((0,\tau],X_\gamma)$ for any $\gamma\in [\frac{1}{2},1)$.

\end{lem}

\textit{Proof :} Note that the initial data is singular with respect to the 
nonlinearity since the latter needs a first derivative but $W_0 \in X_0 $. Although the 
argument is relatively well known, we give the proof for completeness. We shall 
adapt an argument given for instance in \cite[theorem 51.25, p.495]{QS}.\\
We define $E=E_1 \cap E_2$, where 
$E_1=C([0;\tau ]; X_0)$,\\
$ E_2=\{ W \in L^\infty_{loc}((0;\tau], X_{\frac{1}{2}}), \; \underset{ t\in (0;\tau] }
{\sup} 
\sqrt{t} \;\| W(t) \| _ {C^1} < \infty  \}$ \\ 
and $\tau$ will be made precise later.
For $W \in E$, we define its norm : 
$$\|W\| _E = \max \left[ \underset{t\in [0;\tau]}{\sup} 
\;\|W(t)\|_\infty ,\underset{ t\in (0;\tau] }{\sup} \sqrt{t} \; \| W(t) \| _{C^1} 
\right]$$
and  for $K\geq 0$ to be made precise later, we set $E_K=\{W\in E,\; \|W\|_E\leq K\}$.\\
$E_K$ equipped with the metric induced by $\|\;\|_E$ is a complete space.\\
We now define $\Phi :  E_K \longrightarrow E$ by
$$\Phi(W)(t)=S(t)W_0+\int_0^t S(t-s) F_\epsilon(W(s))ds  $$
where $$F_\epsilon(W)=N^2 (m+W) f_\epsilon(m+W+\frac{y.\nabla 
W}{N})$$
For the proof that $\Phi(W)\in C([0;\tau ]; X_0)\bigcap C((0,\tau],X_\gamma)$ for any 
$\gamma\in [\frac{1}{2},1) $ when $W\in E$, we refer to \cite{QS}, p.496 since the proof 
is similar.\\
Next, by properties of analytics semigroups and due to $0<q<2$, we get that for $t\in( 0,
\tau]$ and $W\in E_K$,
$$ \| \Phi(W(t))\|_\infty \leq C_D \| W_0\|_\infty+\frac{\tau^{1-\frac{q}{2}}}{1-\frac{q}
{2}} C_DN^2(m+K)(m\sqrt{\tau}+K)^q$$

$$ \sqrt{t}\| \Phi(W(t))\|_{C^1} \leq C_D \| W_0\|_\infty+\tau^{1-\frac{q}{2}}
C_D N^2(m+K)(m\sqrt{\tau}+K)^qI\left(\frac{1}{2},\frac{q}{2}\right)$$

It is now obvious that $\Phi$ sends $E_K$ into $E_K$ provided that $K\geq 2C_D\| W_0\|
_\infty$ and $\tau$ is small enough. \\
Let $(W_1,W_2)\in (E_K)^2$. We have 
\begin{eqnarray}
 F_\epsilon(W_1)- F_\epsilon(W_2)=&N^2f_\epsilon(m+W_1+\frac{y.\nabla W_1}{N})
[W_1-W_2]&\nonumber\\
\nonumber
 +N^2(W_2+m)&\left[f_\epsilon(m+W_1+\frac{y.\nabla W_1}{N})-
f_\epsilon(m+W_2+\frac{y.\nabla W_2}{N})\right]&
\end{eqnarray}

Now, since $f_\epsilon\in C^1(\mathbb{R})$ and $|f_\epsilon'|\leq L_\epsilon$, we see 
that for any 
$s\in(0,\tau]$,
\begin{eqnarray}
\nonumber
& \|F_\epsilon(W_1(s))- F_\epsilon(W_2(s))\|_\infty\leq& 
N^2\frac{(m\sqrt{s}+K)^q}{s ^{\frac{q}{2}}}\|(W_1-W_2)(s)\|_\infty\\
&&+N^2(m+K)L_\epsilon \|
(W_1-W_2)(s)\|_{C^1}\nonumber \\ 
&&\leq \beta_1(s) \|
W_1-W_2\|_E\nonumber
\end{eqnarray}
where $\beta_1(s)= N^2\left[\frac{(m\sqrt{s}+K)^q}{s ^{\frac{q}
{2}}}+(m+K)\frac{L_\epsilon}{\sqrt{s}}\right]$.

Let $t\in (0,\tau]$. \\
Since $\Phi(W_1)(t)-\Phi(W_2)(t)=\int_0^t S(t-s)[F_\epsilon(W_1(s))- 
F_\epsilon(W_2(s))]ds$, we have
$$\|\Phi(W_1)(t)-\Phi(W_2)(t)\|_\infty \leq  \beta_2 \|W_1-W_2\|_E$$
$$\sqrt{t}\|\Phi(W_1)(t)-\Phi(W_2)(t)\|_{C^1} \leq \beta_3 \|W_1-W_2\|_E$$
where 
$$\beta_2=C_D N^2\left[ \frac{\tau^{1-\frac{q}{2}}}{1-\frac{q}
{2}}  (m\sqrt{\tau}+K)^q
+2\sqrt{\tau} (m+K)L_\epsilon\right]$$
and $$\beta_3=C_D N^2\left[ \tau^{1-\frac{q}{2}}  
(m\sqrt{\tau}+K)^qI(\frac{1}{2},\frac{q}{2})
+\sqrt{\tau} (m+K)L_\epsilon I(\frac{1}{2},\frac{1}{2})\right].$$

So, since $0<q<2$, $\Phi$ is a contraction for $\tau$ small enough. Hence, there exists a 
fixed 
point of $\Phi$, that is to say a mild solution. \\
The uniqueness of the mild solution comes from the uniqueness of the fixed point given by 
the 
contraction mapping theorem.

\begin{rmq} \quad 
\label{rmq_alpha_C1}
If $W_0\in X_{\frac{1}{2}}=\{W\in C^1(\overline{B}),\; W|_{\partial B}=0\}$, then a 
slight modification of the proof shows 
that $W^\epsilon\in C([0,\tau],C^1(\overline{B}))$. Indeed,
we just have to replace the space $E$ in the proof by $E=C([0,\tau],X_{\frac{1}{2}})$. 
Or, we also 
can  refer to \cite[theorem 51.7, p.470]{QS}.
This remark will be helpful later for a density argument.

\end{rmq}

\begin{lem}Let $\epsilon>0$ and $W_0 \in X_0 $.
\begin{itemize}
 \item [i)]
 There exists $T^*_\epsilon=T^*_\epsilon(W_0)>0$ and a unique maximal 
\begin{equation}
\label{alpha_app_1}
W^\epsilon\in C([0,T^*_\epsilon),X_0)\bigcap C^1((0,T^*_\epsilon),X_0)\bigcap 
C((0,T^*_\epsilon),X_1) 
\end{equation}
 such that $ W^\epsilon(0)=W_0$ and for all $t\in (0,T^*_\epsilon)$,
\begin{equation} 
\label{alpha_app_2}
 \frac{d}{dt}W^\epsilon(t)=\Delta W^\epsilon(t) + N^2 (m+W^\epsilon(t)) 
f_\epsilon\left[m+W^\epsilon(t)+\frac{y.\nabla 
W^\epsilon(t)}{N}\right] 
\end{equation}
\begin{equation}
\label{alpha_app_3}
 \underset{t\in(0,T]}{\sup} \sqrt{t}\; \|W^\epsilon(t)\|_{C^1(\overline{B})}<\infty 
 \mbox{ for all } T\in(0,T^*_\epsilon)
\end{equation}
\item [ii)]
Moreover, if $0<t_0<T<T^*_\epsilon$, then for any $\gamma \in(0,1)$,  
$$W^\epsilon\in C^{\frac{1+\gamma}{2},2+\gamma} ( [t_0,T]\times \overline{B})$$
\item [iii)]In particular,  $W^\epsilon\in C ( [0,T^*_\epsilon)\times 
\overline{B})\bigcap 
C^{1,2} ( (0,T^*_\epsilon)\times \overline{B})$ is the unique maximal classical solution 
of the following problem  :
\begin{eqnarray}
\label{alpha_app-1}
  W_t=\Delta W + N^2 (m+W) f_\epsilon[m+W+\frac{y.\nabla 
W}{N}]  &\mbox{ on }(0,T^*_\epsilon)\times \overline{B}&\\
\label{alpha_app-2}
 W=0  &\mbox{ on } [0,T^*_\epsilon)\times \partial B&\\
\label{alpha_app-3}
W(0)=W_0&&
\end{eqnarray}
Moreover, $W^\epsilon$ satisfies 
\begin{equation}
\label{alpha_app-4}
 \underset{t\in(0,T]}{\sup} \sqrt{t}\; \|W^\epsilon(t)\|_{C^1(\overline{B})}<\infty 
 \mbox{ for all } T\in(0,T^*_\epsilon)
\end{equation}
and we have the following blow-up alternative :
$$ T^*_\epsilon=+\infty \quad \mbox{ or } \quad  \underset{t \rightarrow T^*_\epsilon}
{\lim} \|
W^\epsilon(t) \|_{\infty}=+\infty$$
\end{itemize}
\end{lem}

Proof : 
i) In the proof of Lemma \ref{lem_w_app}, we notice that for fixed $\epsilon>0$, the 
minimal existence time $\tau=\tau(\epsilon)$  is uniform for all 
$W_0\in X_0$ such that $\|W_0 \|_\infty\leq r$, where $r>0$.
Then a standard argument shows that there exists a unique maximal mild solution 
$W^\epsilon$ with existence time $T_\epsilon^*>0$ of problem $(\ref{pb_auxiliaire_1})
(\ref{pb_auxiliaire_2})(\ref{pb_auxiliaire_3})$.
It also gives the following blow-up alternative :
$$ T^*_\epsilon=+\infty \quad \mbox{ or }\quad  \underset{t \rightarrow T^*_\epsilon}
{\lim} \|
W^\epsilon(t) \|_{\infty}=+\infty$$
For reference, see for instance \cite[Proposition 16.1, p. 87-88]{QS}.\\
Clearly, $W^\epsilon$ satisfies $(\ref{alpha_app_3})$.
Let us show that $W^\epsilon$ satisfies $(\ref{alpha_app_1})(\ref{alpha_app_2})$.\\ 
Let $t_0\in (0,T^*_\epsilon)$, $T\in (0,T^*_\epsilon-t_0)$ and $t\in [0,T]$. Then,
$$ W^\epsilon (t_0+t)=S(t)W^\epsilon(t_0)+\int_0^t S(t-s)
F_\epsilon(W^\epsilon(t_0+s))ds $$
Since $W^\epsilon\in C((0,T_\epsilon^*),X_{\frac{1}{2}})$, we have $\underset{0\leq s\leq 
T}{\sup}\|W^\epsilon(t_0+s)\|_{C^1}  <\infty$. Then, 
$$F_\epsilon(W^\epsilon(t_0+\cdot))\in L^\infty((0,T),X_0)$$
 We now apply \cite[proposition 4.2.1, p.129]{Lu} to get $W^\epsilon(t_0+\cdot)\in 
 C^{\frac{1}{2}}([0,T],X_{\frac{1}{2}})$. 
Then, $$F_\epsilon(W^\epsilon(t_0+\cdot))\in C^{\frac{1}{2}}([0,T],X_0)$$
We eventually apply \cite[theorem 3.2, p.111]{Pa} to conclude that 
$W^\epsilon$ satisfies $(\ref{alpha_app_1})(\ref{alpha_app_2})
(\ref{alpha_app_3})$ on any 
segment $[t_0,T]\subset (0,T^*_\epsilon)$, hence on $(0,T^*_\epsilon)$.\\
Conversely, since a solution of $(\ref{alpha_app_1})(\ref{alpha_app_2}) 
(\ref{alpha_app_3})$ is a mild
solution, this proves the maximality and the uniqueness.\\
ii)  Let $t_0\in (0,T^*_\epsilon)$, $T\in (0,T^*_\epsilon-t_0)$ and $t\in [0,T]$.
Since in particular $F_\epsilon(W^\epsilon(t_0+\cdot))\in C([0,T],X_0)$, then
$F_\epsilon(W^\epsilon(t_0+\cdot))\in L^p( [t_0,T]\times \overline{B})$ for any $p\geq 
1$. 
Hence, since $W^\epsilon$ satisfies (\ref{alpha_app_2}), by 
interior boundary $L^p$-estimates, we obtain that $W^\epsilon\in W^{1,2}_p ( 
[t_0,T]\times 
\overline{B})$ for any $1\leq p<\infty $. Hence, by Sobolev's  embedding theorem, (see 
for instance \cite[p.26]{Hu}) we see that
$$W^\epsilon\in C^{\frac{1+\gamma}{2},1+\gamma} ( [t_0,T]\times \overline{B}) \mbox{ for 
any 
$\gamma \in(0,1)$ }$$
Eventually, since
$F_\epsilon(W^\epsilon(t_0+\cdot))\in C^{\frac{\gamma}{2},\gamma} ( [t_0,T]\times 
\overline{B})$ then by Schauder interior-boundary parabolic estimates,
$$W^\epsilon\in C^{1+\frac{\gamma}{2},2+\gamma} ( [t_0,T]\times \overline{B}) \mbox{ for 
any 
$\gamma \in(0,1)$}$$
iii) Allmost all is obvious now. Since a classical solution is mild and a mild solution 
is 
classical as seen in i) and ii), then $W^\epsilon$ is also maximal in the sense of the 
classical solutions of $(\ref{alpha_app-1})(\ref{alpha_app-2})$ and 
$(\ref{alpha_app-3})$. \\
The uniqueness of the maximal classical solution comes from the uniqueness of the maximal 
mild solution.\\

\underline{\textit{Proof of theorem \ref{thm_existence_w_epsilon}} :} i)ii)v) The 
correspondence between the solutions of  problem $(\ref{alpha_app-1})
(\ref{alpha_app-2})(\ref{alpha_app-3})$ and problem 
$(tPDE_m^\epsilon)$ is given by $w^\epsilon=W^\epsilon+m$. The 
previous lemma then gives the result. Note that the  existence time is of course the same 
for 
both problems.\\
iii) Let us set $L=\max(K,m)$
 and for $(t,x)\in [0,\frac{1}{qL^qN^2})\times \overline{B}$,
 $$\overline{w}(t,x)=\frac{L}{(1-qL^qN^2t)^{\frac{1}{q}}}.$$  
 Obviously,
$\overline{w}(t)|_{\partial B}\geq L\geq m$ for all $t\geq 0$ and $\overline{w}(0)=L\geq 
w_0$. \\
Moreover, $\overline{w}_t=N^2\overline{w}^{q+1}\geq 
\Delta\overline{w}+N^2\overline{w}f_\epsilon(\overline{w}) $ since for all $x\in 
\mathbb{R}$, $|
f_\epsilon(x)|\leq |x|^q$.\\
Then $\overline{w}$ is a supersolution for problem $(tPDE_m^\epsilon)$, 
so if $0\leq t<\min(T^*_\epsilon,\frac{1}{2qL^qN^2})$, then 
$0\leq w^\epsilon(t)\leq \overline{w}(t)\leq 2^{\frac{1}{q}}L$. \\
We set $\tau=\frac{1}{2qN^2L^q} $ and $C=2^{\frac{1}{q}}L $.
By blow-up alternative  ii), we get
$$T^*_\epsilon\geq \tau \mbox{ \qquad and \qquad }\underset{t \in [0,\tau]}{\sup}\|
w^\epsilon(t)\|_{\infty, \overline{B}}\leq  
C$$
Note that $\tau$ and $C$ depend on $K$, but is independent of $\epsilon$.\\
iv) Noting $W_0=w_0-m$, then for $t\in [0,\tau]$,
 $$ w^\epsilon(t)=m+S(t)W_0+\int_0^t S(t-s)N^2w^\epsilon 
 f_\epsilon\left(w^\epsilon+\frac{x.\nabla 
 w^\epsilon}{N}\right)ds$$
 so
 $$ \|w^\epsilon(t)\|_{C^1}\leq m+ \frac{C_D}{\sqrt{t}}(C+m) +N^2\int_0^t \frac{C_D}
 {\sqrt{t-s}}C 
 \|w^\epsilon(s)\|_{C^1}^q ds $$
 Setting $h(t)=\underset{s\in (0,t]}{\sup} \sqrt{s}\|w^\epsilon(s)\|_{C^1}$, we have 
 $h(t)<\infty$ by (\ref{maj_w_eps}) and 
 $$ \sqrt{t}\|w^\epsilon(t)\|_{C^1}\leq m \sqrt{\tau}+ C_D(C+m) +N^2CC_D\sqrt{t}\int_0^t 
 \frac{1}
 {s^{\frac{q}{2}}\sqrt{t-s}} h(s)^q ds $$
 $$ \sqrt{t}\|w^\epsilon(t)\|_{C^1}\leq m \sqrt{\tau}+ C_D(m+C) +N^2CC_D I\left(\frac{1}
 {2},\frac{q}{2}\right) 
 t^{1-\frac{q}{2}} h(t)^q  $$
 Let $T\in(0,\tau]$. Then,
 \begin{equation}
 \label{h(T)}
  h(T)\leq m \sqrt{\tau}+ C_D(m+C) +N^2CC_D I\left(\frac{1}{2},\frac{q}{2}\right) 
  T^{1-\frac{q}{2}} h(T)^q 
   \end{equation}
Setting $A=m \sqrt{\tau}+ C_D(m+C)$ and $B=N^2CC_D I \left( \frac{1}{2},\frac{q}{2} 
\right)  2^q$, assume that 
there exists $T\in [0,\tau]$ such that $h(T)=2A$. Then,
$$A^{1-q}\leq B T^{1-\frac{q}{2}} \mbox{ which implies } T\geq \left(\frac{A^{1-q}}
{B}\right)^
{\frac{1}{1-q}}$$
Let us set $\tau'=\min\left(\tau,\frac{1}{2}\left(\frac{A^{1-q}}{B}\right)^{\frac{1}{1-
q}} \right)$. \\
Since $h\geq 0$ is nondecreasing, $h_0=\underset{t\rightarrow 0^+}{\lim} h(t)$ exists and  
$h_0\leq A$ by (\ref{h(T)}). So by continuity of $h$ on $(0,\tau']$, 
$ h(t)\leq 2A \mbox{ for all }t\in(0,\tau'] $, that is to say :
$$ \|w^\epsilon(t)\|_{C^1}\leq \frac{2A}{\sqrt{t}} \mbox{ for all }t\in(0,
\tau']$$ 
where $A$ and $\tau'$ only depend on $K$.\\
vi) $0$ is a subsolution of problem 
$(tPDE_m^\epsilon)$.  Then by  comparison principle $w^\epsilon\geq 0$.
The strong maximum principle implies that $w^\epsilon >0$ on $(0,T^*_{\epsilon'})\times 
\overline{B}$ (see  \cite[theorem 5, p.39]{Friedman}). \\
vii) This fact is a consequence of remark \ref{rmq_alpha_C1}.

\subsection{Solutions of problem $(PDE^\epsilon_m)$}

Using the connection with problem $(tPDE_m^\epsilon)$ through the transformation 
$\theta_0$, we shall now provide the proof of Theorem \ref{thm_existence_u_epsilon}. \\

\textit{Proof : } i)
The uniqueness of the classical maximal solution for problem $(PDE_m^\epsilon)$ comes 
from the comparison principle for this problem. \\
We shall now exhibit a classical solution of problem $(PDE_m^\epsilon)$ satisfying 
$(\ref{maj_u_eps})$ and will prove in i)bis) that it is maximal.\\
Let us set $w_0=\theta_0(u_0)$. \\
Remind that $w^\epsilon\in C([0,T^*_\epsilon)\times 
\overline{B})\bigcap C^{1,2}((0,T^*_\epsilon)\times 
\overline{B})$.\\
Remark that a classical solution of $(tPDE_m^\epsilon)$ composed with a rotation 
is still a classical solution. Then by uniqueness, since $w_0$ is radial, so is 
$w^\epsilon(t)$ for all $t\in [0,T^*_\epsilon)$.\\
Hence, $w^\epsilon(t,y)=\tilde{w}^\epsilon(t,\|y\|)$ for all $(t,y)\in 
[0,T^*_\epsilon)\times 
\overline{B}$.
Let us define :
\begin{equation}
 \label{def_u_epsilon}
u^\epsilon(t,x)=x \;\tilde{w}^\epsilon\left(\frac{t}{N^2},x^{\frac{1}{N}}\right)\mbox{ 
for }(t,x)\in 
 [0,N^2 T_\epsilon^*)\times [0,1]
\end{equation}

Since
\begin{equation}
 \label{derivee_en_t_u_epsilon}
(u^\epsilon)_t(t,x)=\frac{x}{N^2} \; (\tilde{w}^\epsilon)_{t}
\left(\frac{t}{N^2},x^{\frac{1}{N}}\right)
\end{equation}
\begin{eqnarray}
 \label{derivee_u_epsilon}
(u^\epsilon)_x(t,x)&=&\left[\tilde{w}^\epsilon+\frac{r(\tilde{w}^\epsilon)_r}{N}\right] 
\left(\frac{t}{N^2},x^{\frac{1}{N}}\right)\nonumber\\
&=&\left[\tilde{w}^\epsilon+\frac{y.\nabla \tilde{w}^\epsilon)}{N}\right] 
\left(\frac{t}{N^2},x^{\frac{1}{N}}\right)
\end{eqnarray}
\begin{eqnarray}
 \label{derivee_seconde_u_epsilon}
x^{2-\frac{2}{N}}(u^\epsilon)_{xx}(t,x)&=&\frac{x}
{N^2}\left[(\tilde{w}^\epsilon)_{rr}+\frac{N+1}{r}(\tilde{w}^\epsilon)_{r}\right] 
\left(\frac{t}{N^2},x^{\frac{1}{N}}\right)\nonumber \\
&=&\frac{x}{N^2}\; \Delta w^\epsilon \left(\frac{t}{N^2},x^{\frac{1}{N}}\right)
\end{eqnarray}
 it  is easy to check that $$u^\epsilon\in C([0,N^2T^*_\epsilon)\times 
 [0,1])\bigcap C^1((0,N^2T^*_\epsilon)\times [0,1])\bigcap
 C^{1,2}((0,N^2T^*_\epsilon)\times (0,1])$$ is a classical solution of  problem 
 $(PDE_m^\epsilon)$ on 
$[0,N^2 T^*_\epsilon )$.\\
Special attention has to be paid to the fact that $u^\epsilon$ is $C^1$ up to $x=0$ but 
this is 
clear because of (\ref{def_u_epsilon}) and (\ref{derivee_u_epsilon}).\\
Since $\|w_0\|_{\infty,\overline{B}}=\mathcal{N}[u_0]\leq K$, formula 
(\ref{derivee_u_epsilon}) and theorem \ref{thm_existence_w_epsilon} iii)iv) imply that 
there exist $\tau'=\tau'(K)\in (0,1]$ and $C=C(K)>0$ independent of $\epsilon$ such that 
$\underset{t\in(0,\tau']}{\sup}  \|w^\epsilon(t)\|_{\infty,
\overline{B}}+\underset{t\in(0,\tau']}{\sup} \sqrt{t}\; \|w^\epsilon(t)\|_{C^1,
\overline{B}}\leq C$. Then,  

\begin{equation}
\label{condition_t=0} 
 \underset{t\in(0,\tau']}{\sup} \sqrt{t}\; \|u^\epsilon(t)\|_{C^1([0,1])}\leq C 
 \end{equation}
It is also clear from formula $(\ref{def_u_epsilon} )$ that 
$T_{max}^\epsilon \geq N^2 T_\epsilon^*$.\\
ii)
From theorem \ref{thm_existence_w_epsilon} iv), $T^*_\epsilon\geq \tau'$, then 
$u^\epsilon$  is at least defined on $[0,N^2\tau']$ and can be extended to a maximal 
solution. This minimal existence time $\tau=N^2 \tau'$ only depends on $K$. \\
Moreover, by formula (\ref{def_u_epsilon}),
$$\underset{t\in [0,\tau]}{\sup} \mathcal{N}[u^\epsilon(t)]= \underset{(t,y)\in [0,
\tau']\times \overline{B}\backslash \{0\} }{\sup} w^\epsilon(t,y)=\underset{t\in [0,
\tau'] }{\sup} \|w^\epsilon(t)\|_{\infty,\overline{B}}\leq C$$ \\
i)bis) If $T_\epsilon^*= \infty $, then formula $(\ref{def_u_epsilon})$ gives a global 
solution $u^\epsilon$ then $T_{max}^\epsilon=\infty$.\\
Suppose $T_\epsilon^*< \infty $.\\
Assume that $T_{max}^\epsilon> N^2 T_\epsilon^*$ with maybe $T_{max}^\epsilon =\infty$.
Then, there exists in particular a classical solution $u^\epsilon$ of $(PDE_m^\epsilon)$ 
on $[0,N^2 T_\epsilon^*]$. By uniqueness, on $[0,N^2 T_\epsilon^*)$, 
$u^\epsilon$ coincides with the solution given by 
$(\ref{def_u_epsilon})$.
By blow-up alternative for $w^\epsilon$, 
$\underset{t\rightarrow T^*_\epsilon}{\lim} 
\|w^\epsilon(t)\|_{\infty,\overline{B}} =\infty$ thus 
$\underset{t\rightarrow N^2T^*_\epsilon}{\lim}  
\mathcal{N}[u^\epsilon(t)] =\infty$.\\
But, since (\ref{majoration_u_epsilon_tau'}) and 
$u^\epsilon\in C([\tau',N^2 T^*_\epsilon], C^1( [0,1]))$ 
then 
$$\underset{t\in [0,N^2T^*_\epsilon]}{\sup} 
\mathcal{N}[u^\epsilon(t)] <\infty$$
 which provides a contradiction. Whence i).\\
 Moreover, this proves that the solution 
$u^\epsilon$  is actually maximal.\\
iii) The blow-up alternative for problem $(PDE_m^\epsilon)$ 
follows directly from $i)bis)$ and from the blow-up alternative 
for problem $(tPDE_m^\epsilon)$.\\
iv) This point needs some work that will be done in the next lemma.\\
v) This follows  from Theorem \ref{thm_existence_w_epsilon} iv) and formulas 
(\ref{def_u_epsilon})(\ref{derivee_en_t_u_epsilon})(\ref{derivee_u_epsilon})
(\ref{derivee_seconde_u_epsilon}).\\
vi) This follows  from Lemma \ref{propriete_theta})i), Theorem 
\ref{thm_existence_w_epsilon} vii) and formula (\ref{derivee_u_epsilon}).\\

The next lemma, whose proof is rather technical, is very important since it shows that $(u^\epsilon)_x>0$ on 
$(0,T_{max}^\epsilon)\times [0,1]$, which will imply later that solutions of 
$(PDE_m)$ at time $t$ are nondecreasing. Moreover, this fact is essential in 
\cite{Montaru} in 
order to prove 
that some functional $\mathcal{F}_\epsilon$ is a strict Lyapunov functional 
for the dynamiacl system induced by problem $(PDE_m^\epsilon)$.\\

\begin{lem}
\label{lem_derivee_u_epsilon}
 Let $\epsilon>0$, $u_0\in Y_m$ and $w_0=\theta_0(u_0)$. \\Let us set 
 $T_{max}^\epsilon=T^\epsilon_{max}(u_0)$ and $T_\epsilon^*=T_\epsilon^*(w_0)$.
\begin{itemize}
\item[i)] $0\leq u^\epsilon \leq m$ on $[0,T_{max}^\epsilon)\times [0,1]$.
\item[ii)] $u^\epsilon\in C^{1,3}((0,T_{max}^\epsilon)\times (0,1])$
and 
$u^\epsilon\in C^2((0,T_{max}^\epsilon)\times (0,1])$ (not optimal).
 \item [iii)] For all $(t,x)\in (0,T_{max}^\epsilon)\times [0,1]$, 
 $(u^\epsilon)_x(t,x)>0$.
 \item[iv)]  
 $w^\epsilon+\frac{y.\nabla w^\epsilon}{N}>0$ for any $(t,y)\in (0,T_\epsilon^*)\times 
 \overline{B}$.
\end{itemize}
\end{lem}

Proof : 
$i)$ $0$ and $m$ are respectively sub- and supersolution for problem 
$(PDE_m^\epsilon) $ satisfied by $u^\epsilon$. Whence 
the result by comparison principle.\\

ii) Let $[t_0,T]\subset (0,T_{max}^\epsilon)$.\\
Let $w_0=\theta_0(u_0)$. We set  $t'_0=\frac{t_0}{N^2}$ and $T'=\frac{T}{N^2}$. 
We now refer to \cite[p.72, Theorem 10]{Friedman} and apply it to $D=(t_0',T')\times 
\overline{B}$. \\
We recall that 
$w^\epsilon$ satisfies on $D$
\begin{equation}
\label{equation_w}
w_t=\Delta w+c\;w 
\end{equation}
with $c=N^2\, f_\epsilon(w+\frac{y.\nabla w}{N})$.
Let $\gamma\in(0,1)$. 
$\nabla c$ is H\"older continuous with exponent $\gamma$ in $D$  because $w^\epsilon\in 
C^{1+\frac{\gamma}{2},2+\gamma}([t_0',T']\times \overline{B})$  and 
$f_\epsilon'$ is Lipschitz continuous 
on compact sets of $\mathbb{R}$. Then $\partial_t \nabla w^\epsilon$ and
$\partial_\alpha w^\epsilon$ 
 are 
H\"older continuous with exponent $\gamma$ in $D$ for any multi-index
$|\alpha|\leq 3$. 
Thus, \\
$w\in C^{1,3}([t_0',T']\times \overline{B})$, so $u^\epsilon\in C^{1,3}([t_0,T]\times 
(0,1])$ by formula (\ref{def_u_epsilon}).\\
We apply the same theorem again : $\partial_\alpha c$ is 
H\"older continuous with exponent $\gamma$ for any $|\alpha|\leq 2$ since $f_\epsilon ''$ 
is Lipschitz 
continuous on compact sets of $\mathbb{R}$. So, $\partial_t\partial_\alpha w$ is 
H\"older continuous with exponent $\gamma$ for any  $|\alpha|\leq 2$.
Then, $c_t$ and $\partial_t\Delta w$ are continuous so by (\ref{equation_w}),
$w_{tt}$ is continuous.\\
By (\ref{equation_w}) again, it is clear that $\partial_\alpha \partial_t w$ is 
continuous for $|\alpha|\leq 1$ hence $w\in C^2([t_0',T']\times \overline{B})$. 
It follows from formula (\ref{def_u_epsilon}) that 
$u^\epsilon\in C^2([t_0,T]\times (0,1])$.\\
In particular, $(u^\epsilon)_{t,x}=(u^\epsilon)_{x,t}$.

$iii)$ Let $T\in(0,T_{max}^\epsilon)$.\\
We prove the result in two steps.\\
 \underline{First step :}
 We now show that $v^\epsilon:=u^\epsilon_x \geq 0$ on $[0,T]\times 
 [0,1]$. \\
 We divide the proof in three parts.
 
 \begin{itemize}
 \item 
  \underline{First part :} We show the result for any $u_0\in Y_m^{1,\gamma}$ where 
  $\gamma>\frac{1}
 {N}$.\\
Since $u^\epsilon$ satisfies on $(0,T]\times (0,1]$
\begin{equation}
\label{equation_u_epsilon}
u^\epsilon_t=x^{2-q}u^\epsilon_{xx}+u^\epsilon f_\epsilon(u^\epsilon_x)
\end{equation}
and thanks to ii), 
we can now differentiate this equation with respect to $x$. We denote 
$$b=\left[(2-\frac{2}{N})x^{1-\frac{2}{N}}+u^\epsilon f_\epsilon'(v)\right]$$
  and obtain the partial 
differential equation satisfied by $v^\epsilon$ :
\begin{eqnarray}
\label{v_app_1}
 v_t=x^{2- \frac{2}{N}}v_{xx}+b\;v_x+f_\epsilon(v)v &\mbox{ on }&
(0,T)\times(0,1)\\
\label{v_app_2}
 v(0,\cdotp)=(u_0)'&&\\
 \label{v_app_3}
v(t,0)=u^\epsilon_x(t,0) &\mbox{   for  }& t\in(0,T] \\
\label{v_app_4}
v(t,1)=u^\epsilon_x(t,1) &\mbox{        for } &t\in(0,T] 
\end{eqnarray}

By Theorem \ref{thm_existence_u_epsilon} vii), we know that $u^\epsilon\in C([0,T], C^1( 
[0,1]))$, then $v^\epsilon \in C([0,T]\times [0,1])$
and $v^\epsilon$ reaches its minimum on $[0,T]\times [0,1]$.\\ 
From $i)$ follows that 
$u^\epsilon_x(t,0)\geq 0$ and $u^\epsilon_x(t,1)\geq 0$ for all $t\in(0,T]$.
Then, from $(\ref{v_app_2}),(\ref{v_app_3})$ and $(\ref{v_app_4})$, $v^\epsilon\geq 0$ on 
the parabolic boundary of $[0,T]\times [0,1]$. From $(\ref{v_app_1})$, we see that 
$v^\epsilon$ cannot reach a negative minimum 
in
$(0,T]\times (0,1)$ since for all $x\neq 0$, $x\,f_\epsilon(x)>0$. So 
$v^\epsilon \geq 0$ on $[0,T]\times [0,1]$.

\item
 \underline{Second part :} We show that if $u_0\in Y_m$, there exists $\tau>0$ 
 such that for all
 $t\in [0,\tau]$,  $u^\epsilon(t)$ is non decreasing on $[0,1]$.\\
 Let $u_0\in Y_m$. 
From Lemma \ref{lem_density}, there exists a sequence 
 $(u_n)_{n\geq 1}$  of $ Y_m^{1,1}$  such that 
 $\|u_n-u_0\|_{\infty,[0,1]} 
 \underset{n\rightarrow \infty}{\longrightarrow} 0$ and
$\mathcal{N}[u_n]\leq \mathcal{N}[u_0]$.\\
By Theorem \ref{thm_existence_u_epsilon} ii),
  there exists a common small existence time $\tau>0$ for all solutions $(u_n^\epsilon)_{ 
  n\geq 0}$ of problem $(PDE_m^\epsilon)$ with initial condition $u_n$.
 From first part, we know that for all
 $t\in [0,\tau]$ $u_n^\epsilon(t)$ is a nondecreasing function since $u_n\in Y_m^{1,1}$. 
 To prove the result, it is sufficient to show that $\|u_n^\epsilon-u^\epsilon\|_{\infty,
 [0,1]\times [0, \tau]} 
 \underset{n\rightarrow \infty}{\longrightarrow} 0$.\\
 Let $\eta>0$. By $(\ref{maj_u_eps})$, there exists $C>0$
 such that for all $t\in [0,\tau]$, $\|(u^\epsilon(t))_x\|_\infty \leq 
 \frac{C}{\sqrt{t}}$. So we can choose
  $\eta'>0$ such that $$\eta' e^{\int_0^\tau [\|(u^\epsilon(t))_x\|_\infty^q+1] \;dt}\leq 
  \eta$$
 Let $n_0\geq 1$ such that for all $n\geq n_0$, $\|u_n-u_0\|_{\infty,[0,1]}\leq \eta' $. 
 Let $n\geq n_0$.\\
 Let us set $$z(t)=[u_n^\epsilon(t)-u^\epsilon(t)]e^{-\int_0^\tau [\|(u^\epsilon(t))_x\|
 _\infty^q+1] \;dt}$$ 
 We see that $z$ satisfies
 \begin{equation}
 \label{equ_aux_z}
  z_t=x^{2-\frac{2}{N}}z_{xx}+b\, z_x+c\, z
 \end{equation}
 
 where $b=u_n^\epsilon \, \dfrac{f_\epsilon((u_n^\epsilon)_x)-f_\epsilon((u^\epsilon)_x)}
 {(u_n^\epsilon)_x-(u^\epsilon)_x}$ if 
 $(u_n^\epsilon)_x\neq (u^\epsilon)_x$ and $0$ else,\\
  $c=[f_\epsilon((u^\epsilon)_x)- \|(u^\epsilon)_x\|_\infty^q-1]<0$.\\
  Since $z\in C([0,\tau]\times [0,1])$, $z$ reaches its maximum and its minimum. \\
  Assume that this maximum is greater than $\eta'$. Since $z=0$ for $x=0$ and $x=1$ and 
  $z\leq \eta'$ for $t=0$, it can be reached only in 
  $(0,\tau]\times (0,1)$ but this is impossible because $c<0$ and (\ref{equ_aux_z}). 
  We make the similar reasoning for the minimum.
   Hence,
  $|z|\leq \eta'$ on  $[0,\tau]\times [0,1]$. \\Eventually,
  $\|u_n^\epsilon-u^\epsilon\|_{\infty,[0,1]\times [0, \tau]} \leq \eta'
   e^{\int_0^\tau [\|(u^\epsilon(t))_x\|_\infty^q+1] \;dt}\leq \eta$ for all $n\geq n_0$. 
   Whence the result.
   
   \item
 \underline{Last part :} Let $u_0\in Y_m$. From the second part, there exists $\tau>0$ 
 such that that for all $t\in[0,\tau]$, $u^\epsilon(t)$ is nondecreasing. Since 
 $u^\epsilon\in C([\tau,T_{max}^\epsilon),C^1([0,1])  )$ and $u_0(\tau)$ is 
 nondecreasing, we can apply the same argument as in the first part to deduce that for 
 all $t\in [\tau,T_{max}^\epsilon)$, $u^\epsilon(t)$ is nondecreasing. That concludes the 
 proof of the second step.
 \end{itemize}

 \underline{Second step :} Let us show that $v^\epsilon>0$ on $(0,T]\times [0,1]$.\\
First, from formula (\ref{derivee_u_epsilon}) and Theorem 
\ref{thm_existence_w_epsilon} vi) follows that  $v^\epsilon(t,0)=(u^\epsilon)_x(t,0)>0 $ 
for 
$t\in(0,T]$.\\
Assume by contradiction that $v^\epsilon$ is zero at some point
 in $(0,T)\times (0,1)$. \\
 Let $z=v^\epsilon e^{-\int_0^t[ \|v^\epsilon(s)\|_\infty^q+1]}ds\geq 0$ by second step. 
 $z$ reaches its minimum and satisfies the following equation :
\begin{equation}
z_t=x^{2-\frac{2}{N}}z_{xx}+[(2-\frac{2}{N})x^{1-\frac{2}{N}}+u^\epsilon 
f_\epsilon'(v^\epsilon)]z_x+[f_\epsilon(v^\epsilon)-
\|v^\epsilon(s)\|_\infty^q-1]z
\end{equation} 
 where $f_\epsilon(v^\epsilon)-
\|v^\epsilon(s)\|_\infty^q-1\leq -1$ on $[0,T]\times [0,1]$.\\
  Then, by the strong minimum principle (\cite{Friedman}, p.39, Theorem 5) 
applied to $z$, 
we deduce that $v^\epsilon=0$ on $(0,T)\times (0,1)$.
Then, by continuity, $v^\epsilon(t,0)=0$ for $t\in(0,T)$ which contradicts the previous 
assertion.\\
Suppose eventually that $v^\epsilon(t,1)=0$ for some $t\in(0,T)$. From $ 
(\ref{equation_u_epsilon})$, we 
deduce
that $(u^\epsilon)_{xx}(t,1)=0$, ie $v^\epsilon_x(t,1)=0$. \\
Since  $f_\epsilon(y)y\geq 0$ for all $y\in \mathbb{R}$, we observe that $v^\epsilon$ 
satisfies :
 \begin{equation}
v_t\geq x^{2-\frac{2}{N}}v_{xx}+[(2-\frac{2}{N})x^{1-\frac{2}{N}}+u^\epsilon 
f_\epsilon'(v)]v_x
\end{equation} 
Since $v^\epsilon>0$ on $(0,T)\times [\frac{1}{2},1)$ and the underlying operator in the 
above 
equation is uniformly parabolic on
$(0,T)\times [\frac{1}{2},1]$, we can apply Hopf's minimum principle (cf. \cite[Theorem 
3, p.170]{PW}) to deduce that 
$v^\epsilon_x(t,1)<0$ what yields a contradiction. In conclusion, $(u^\epsilon)_x>0$ on 
$(0,T]\times [0,1]$ 
for all $T<T_{max}^\epsilon$, whence the result.\\

iv) It is clear from iii) thanks to formula (\ref{derivee_u_epsilon}).\\

We can now deduce the following monotonicity property which will be useful in order to 
find a solution of problem $(PDE_m)$ by letting $\epsilon$ go to zero.

\begin{lem} Let $u_0\in Y_m$ and $w_0=\theta_0(u_0)$.
\label{monotonie_w_epsilon}
\begin{itemize}
\item[i)] If $\epsilon'<\epsilon$, then $T^*_{\epsilon'}\leq T^*_\epsilon$ and 
$w^{\epsilon'}\geq 
w^\epsilon$ on $ [0,T^*_{\epsilon'})\times \overline{B}$.
\item[ii)] If $\epsilon'<\epsilon$, then $T_{max}^{\epsilon'}\leq T_{max}^\epsilon$ and 
$u^{\epsilon'}\geq u^\epsilon$ on $ [0,T_{max}^{\epsilon'})\times [0,1]$.
\end{itemize}
\end{lem}

\textit{Proof :}
i) $w^{\epsilon'}$ is a supersolution for  $(tPDE_m^\epsilon)$
since $w^{\epsilon'}+\frac{y.\nabla w^{\epsilon'}}{N}\geq 0$ for all $(t,y)\in 
[0,T^*_{\epsilon'})\times \overline{B}$ and $f_{\epsilon'}\geq f_\epsilon$ on $[0,
+\infty)$ for $\epsilon' <\epsilon $.
Using the blow-up alternative for problem  $(tPDE_m^{\epsilon})$, we get the result by 
contradiction.\\
ii) It is clear from i) using the relation between $u^\epsilon$ and 
$w$ in Theorem 
$\ref{thm_existence_u_epsilon}$ iii). We could as well use a comparison argument as in 
i).

\begin{rmq} 
 $(w^\epsilon)_{\epsilon\in(0,1)}$ (resp.  $(u^\epsilon)_{\epsilon\in(0,1)}$) is then a 
 nondecreasing family of functions for $\epsilon$ 
decreasing, with an existence time maybe shorter and shorter  but not less than a given 
$\tau>0$ depending on $\|w_0\|_\infty$ (resp. $\mathcal{N}[u_0]$).

\end{rmq}

\subsection{Solutions of problem $(tPDE_m)$ and proof of Theorem \ref{thm_existence_u} }

We shall now prove Theorem \ref{thm_existence_w}, i.e. the local in time wellposedness of 
problem $(tPDE_m)$. \\

The small time existence part is obtained by passing to the limit $\epsilon$ to 0 in 
problem $(tPDE_m^\epsilon)$ via the following lemma :

\begin{lem} \textbf{Local existence of a classical solution for problem $(tPDE_m)$}\\
\label{existence_locale_w} 
Let $w_0=\theta_0(u_0)$ where $u_0\in Y_m$ with $\|w_0\|_{\infty,\overline{B}}\leq K$.\\
There exists $\tau'=\tau'(K)>0$ and $w\in C([0,\tau']\times \overline{B}) \bigcap C^{1,2} 
((0,\tau']\times 
\overline{B})$ such that $w^\epsilon \underset{\epsilon \rightarrow 0}{\longrightarrow}w$ 
in $C([0,
\tau']\times \overline{B})$ and in $ C^{1,2} ((0,\tau']\times \overline{B})$.\\
Moreover, $w$ is the unique classical solution of problem $(tPDE_m)$ on $[0,\tau']$
and satisfies the following condition :
$$ \underset{t\in(0,\tau']}{\sup} \sqrt{t}\; \|w(t)\|_{C^1([0,1])}<\infty $$
\end{lem}

\textit{Proof :} \\
\underline{First step : } From Theorem \ref{thm_existence_w_epsilon} $iii)iv)$,
there exists $\tau'=\tau'(K)>0$ 
and $C=C(K)>0$  both independent of $\epsilon$ such that
$$ \underset{t \in [0,\tau']}{\sup} \|w^\epsilon(t)\|_{\infty,\overline{B}}\leq 
C $$
$$\underset{t \in (0,\tau']}{\sup} \sqrt{t}\|w^\epsilon(t)\|_{C^1(\overline{B})}\leq 
C  $$
Let $t_0\in (0,\tau']$. Recall that $ F_\epsilon(w)=N^2 w 
f_\epsilon(w+\frac{x.\nabla 
w}{N})$.\\
We see that for all $\epsilon>0$ and $t\in [t_0,\tau']$, $\|w^\epsilon(t)\|_{C^1}\leq 
\frac{C}{\sqrt{t_0}}$ where $C$ is independent of $\epsilon$. \\
If $x\geq 0$, $0\leq 
f_\epsilon(x)\leq x^q$, so
 there exists $C'>0$ which depends on $t_0$ but is independent of $\epsilon$ such that
 $$\|F_\epsilon(w^\epsilon)\|_{\infty,[t_0,\tau']\times\overline{B}}\leq C'$$ 
then for any $p\geq 1$, 
$\|F_\epsilon(w^\epsilon)\|_{ L^p([t_0,\tau']\times\overline{B})}\leq C''$
where $C''$ depends on $t_0$ but is independent of  $\epsilon$.\\
We can 
now use the $L^p$ estimates, then Sobolev embedding and eventually interior-boundary 
Schauder 
estimates to obtain that for any $\gamma\in [0,q)$, 
$$\|w^\epsilon\|_{C^{1+\frac{\gamma}{2},2+\gamma}([t_0,\tau']\times\overline{B} )}\leq 
C''' $$
where $C'''$ depnds on $t_0$ but is independent of $\epsilon$ since $f_\epsilon$ is 
H\"older continuous with exponent $q$ on $[0,+\infty)$ and H\"older coefficient less or 
equal to 1.\\ 
We now use a sequence $t_k \underset{k\rightarrow \infty}{\longrightarrow 0}$ and the 
Ascoli's 
theorem for each $k$ and eventually proceed to a diagonal extraction to get a sequence 
$\epsilon_n \underset{n\rightarrow \infty}{\longrightarrow} 0$ such that 
$$ w^{\epsilon_n}  \underset{n\rightarrow \infty}{\longrightarrow}  w$$ in $C^{1,2}([t_k,
\tau']\times\overline{B} )$ for some function $w$, for each $k$. So, 
\begin{equation}
\label{w_loc_1}
w\in C^{1,2}((0,\tau']\times\overline{B} )
\end{equation}
Since by Lemma \ref{monotonie_w_epsilon} i), $w^\epsilon$ is nondecreasing as $\epsilon$ 
decreases to $0$, then $w=\underset{\epsilon 
\rightarrow 0^+}{\lim} \;w^\epsilon$ on $(0,\tau']\times\overline{B} $.
 $w$ is then unique. Hence,
 $$ w^{\epsilon}  \underset{\epsilon \rightarrow 0}{\longrightarrow}  w \mbox{\qquad    
 in }C^{1,2}
 ([t_0,\tau']\times\overline{B} ) \mbox{   for each } t_0\in(0,\tau'] $$
For a fixed $s\in(0,\tau']$,
$w^{\epsilon}(s)  \underset{\epsilon \rightarrow 0}{\longrightarrow}  w(s) \mbox{\;  in 
}C^{1}
(\overline{B} )$, then Lemma \ref{lem_derivee_u_epsilon} iv) implies
\begin{equation}
\label{w_loc_2}
w+\frac{y.\nabla w}{N}\geq 0 \mbox{ \quad on }(0,\tau']\times \overline{B}
\end{equation}

 Moreover, the both following estimates are clear :
 \begin{equation}
\label{w_loc_5}
\|w(t)\|_{\infty,\overline{B}}\leq C \mbox{ for all }t\in(0,
\tau']
\end{equation}
 \begin{equation}
\label{w_loc_3}
\|w(t)\|_{C^1}\leq \frac{C}{\sqrt{t}} \mbox{ for all }t\in(0,
\tau']
\end{equation}
 
\underline{Second step : } Let us show that $w\in C([0,\tau']\times\overline{B})$ and 
that 
 $$ w^{\epsilon}  \underset{\epsilon \rightarrow 0}{\longrightarrow}  w \mbox{\qquad    
 in }C([0,
 \tau']\times\overline{B} )$$
First, remark that from Dini's theorem, the second part is obvious once the first one is 
known 
since $w^\epsilon$ is nondecreasing on the compact set $[0,\tau']\times\overline{B} $ and 
$w^\epsilon$ converges pointwise to the continuous function $w$.\\
Let $t\in(0,\tau']$. Let us set $W_0=w_0-m$. We have
$$ w^\epsilon(t)=m+S(t)W_0+\int_0^t S(t-s)F_\epsilon(w^\epsilon(s)) ds$$
$(w^\epsilon+\frac{y.\nabla w^\epsilon}{N})(t,x)\geq 0$ hence $f_\epsilon( 
w^\epsilon+\frac{y.\nabla w^\epsilon}{N})=( 
w^\epsilon+\frac{y.\nabla w^\epsilon}{N}+\epsilon)^q-\epsilon^q$.\\

Let $s\in(0,t)$. Clearly, $F_\epsilon(w^\epsilon(s)) \underset{\epsilon \rightarrow 0 }
{\longrightarrow} 
N^2w(s)\left(w(s)+\frac{y.\nabla w(s)}{N}\right)^q \mbox{\;  in }C(\overline{B} ) $.\\
By continuous dependence of the heat semi-group on $C_0(\overline{B} )$ with respect to 
the initial data, we have
$$S(t-s) F_\epsilon(w^\epsilon(s)) \underset{\epsilon \rightarrow 0 }{\longrightarrow} 
S(t-s) 
N^2w(s)\left(w(s)+\frac{y.\nabla w(s)}{N}\right)^q \mbox{\;  in }C(\overline{B} ) $$
Moreover, we have a uniform domination for all $\epsilon\in (0,1)$ since
$$ \|S(t-s) F_\epsilon(w^\epsilon(s))\|_{\infty,\overline{B}}
\leq C_D N^2 C (C+\frac{C}{\sqrt{s}})^q$$
and the RHS belongs to $L^1(0,t)$.
Hence, since  
$w^{\epsilon}(t)  \underset{\epsilon \rightarrow 0}{\longrightarrow}  w(t) \mbox{\;  in 
}C(\overline{B} )$, by the Lebesgue's dominated convergence theorem, we obtain :
$$ w(t)=m+S(t)W_0+\int_0^t S(t-s) N^2w(s)\left(w(s)+\frac{x.\nabla w(s)}{N}\right)^q  
ds$$
Then $\|w(t)-w_0\|_{\infty,\overline{B}}\leq \|S(t)W_0-W_0\|_{\infty,
\overline{B}}+
\int_0^t C_D N^2 C (C+\frac{C}{\sqrt{s}})^q$.\\
Hence, by the continuity of the heat semigroup at $t=0$ on  $C_0(\overline{B} )$,
$$w(t)  \underset{t \rightarrow 0}{\longrightarrow}  w_0 \mbox{\;  in }C(\overline{B} )$$ 

We can then deduce
\begin{equation}
\label{w_loc_4}
w\in C([0,\tau'],C(\overline{B}))=C([0,\tau']\times \overline{B})
\end{equation}
\underline{Last step :} Passing to the limit, since 
$ w^{\epsilon}  \underset{\epsilon \rightarrow 0}{\longrightarrow}  w \mbox{ in }C^{1,2}
([t_0,\tau']\times\overline{B} )$ for each $ t_0\in(0,\tau'] $, then 
w satisfies :
$$ w_t=\Delta w+N^2 \,w\,(w+\frac{y.\nabla w}{N})^q \mbox{\quad on \quad} (0,
\tau']\times\overline{B}.$$
Since, moreover, (\ref{w_loc_1})(\ref{w_loc_2})(\ref{w_loc_3})(\ref{w_loc_4}) hold, $w$ 
is thus a classical solution of problem $(tPDE_m)$. The uniqueness
comes from the comparison principle.

\subsection{Proofs of Theorem \ref{thm_existence_w} and Theorem \ref{thm_existence_u}}

\underline{\textit{Proof of Theorem \ref{thm_existence_w} :}} i) and ii) are standard 
since the small existence time depends on $\|w_0\|_{\infty, \overline{B}}$. 
For reference, see  \cite[Proposition 16.1, p. 87-88]{QS} for instance.\\
iii) By Lemma \ref{lem_derivee_u_epsilon} iv), since $f_\epsilon(s)\leq s ^q$ for all 
$s\geq0 $, so $w^\epsilon$ is a subsolution of $(tPDE_m)$ so by 
comparison principle, $$0\leq w^\epsilon\leq w
\mbox{ on } (0,\min(T^*_\epsilon, T^*))$$
 By blow-up alternative for 
classical solutions of $(tPDE_m^\epsilon)$, it is easy to see by contradiction that 
$T^*_\epsilon\geq T^*$. It implies that
$w\geq w^\epsilon>0$ on $(0,T^*)\times \overline{B} $ by Theorem 
\ref{thm_existence_w_epsilon} vi).\\
iv) We use interior-boundary Schauder estimates.\\

\underline{\textit{Proof of Theorem \ref{thm_existence_u} :}} it follows from Theorem 
\ref{thm_existence_w} by exactly  the same way as for passing from Theorem 
\ref{thm_existence_w_epsilon} to Theorem \ref{thm_existence_u_epsilon}. \\
The part vi) will be proved in  subsection \ref{regularite}.

\begin{rmq} We can precisely describe the connection between problems $(PDE_m)$ and 
$(tPDE_m)$. \\
Let $w_0=\theta_0(u_0)$ with $u_0\in Y_m$. Then,
 $$T_{max}(u_0)=N^2 T^*(w_0).$$
Moreover, for all $(t,x)\in [0,T_{max})\times [0,1]$, 
\begin{equation}
\label{def_u}
u(t,x)=x\; \tilde{w}(\frac{t}{N^2},x^{\frac{1}{N}}).
\end{equation}
\end{rmq}

\subsection{Convergence of maximal classical solutions of problem $(PDE_m^\epsilon)$ 
to classical solutions of problem $(PDE_m)$ as $\epsilon$ goes to $0$}

\underline{\textit{Proof of Lemma \ref{lem_convergence_u_epsilon} :}} i) Since $0\leq 
f_\epsilon(s)\leq s^q$ for all $s\geq 0$ and $(u^\epsilon)_x\geq 0$, it is easy 
to check that $u^\epsilon$ is a subsolution for problem $(PDE_m)$ with initial condition 
$u_0$. Hence,
by blow-up alternative, this implies that $T_{max}(u_0)\leq T_{max}^\epsilon(u_0)$.\\
ii)  Let $w_0=\theta_0(u_0)$.\\
We know that $w^\epsilon$  is a subsolution for problem $(tPDE_m)$ with initial condition 
$w_0$ 
thus, setting $T'=\frac{T}{N^2}$, 
$t_0'=\frac{t_0}{N^2}$,
$$\underset{t\in [0,T']}{\sup}\|w^\epsilon(t)\|_{\infty,\overline{B}}\leq 
\underset{t\in [0,T']}{\sup}\|w(t)\|_{\infty,\overline{B}}=:K<\infty.$$
Applying Theorem \ref{thm_existence_w_epsilon} iv), we know that there exists 
$\tau'\in(0,t_0')$ and  $C>0$ 
both depending on $K$ such that for all $w_0\in Z_m$ with $\|w_0\|_{\infty,
\overline{B}}\leq K$ we have
$$\underset{t\in (0,\tau']}{\sup}\sqrt{t}\|w^\epsilon(t)\|_{C^1(\overline{B})}\leq C.$$
So, for all $w_0\in Z_m$ with $\|w_0\|_{\infty,\overline{B}}\leq K$, 
$$\|w^\epsilon(\tau')\|_{C^1(\overline{B})}
\leq \frac{C}{\sqrt{\tau'}}=:C'.$$ where $C'$ depends on $K$ and $t_0$.\\
For $t\in [t_0',T']$, we can use  $w^\epsilon(t-\tau')$ as initial data  to show that
$$\underset{t\in [t_0',T']}{\sup}\|w^\epsilon(t)\|_{C^1(\overline{B})}\leq C'.$$
We can then proceed as in the proof of Lemma \ref{existence_locale_w} and show that
$$w^\epsilon \underset{\epsilon \rightarrow 0}
{\longrightarrow} w \mbox{ in } C^{1,2}([t_0',T']\times \overline{B}).$$
Whence the results thanks to formulas (\ref{def_u_epsilon})(\ref{derivee_en_t_u_epsilon})
(\ref{derivee_u_epsilon})(\ref{derivee_seconde_u_epsilon}) and their equivalent for $u$ 
and  $\tilde{w}$.

\subsection{Regularity of classical solutions of problem $(PDE_m)$}
\label{regularite}
We already know that classical solutions verify $u(t)\in C^1([0,1])$ for all 
$t\in(0,T_{max}(u_0))$ but we can 
actually be more precise, as stated in the next lemma which corresponds exactly to
Theorem \ref{thm_existence_u} vi).

\begin{lem}
Let $u_0\in Y_m$. \\
For all $t\in (0,T_{max}(u_0))$, $u(t)\in Y_m^{1,\frac{2}{N}}$.
\end{lem}

\textit{Proof :} Let $(t,x)\in (0,T_{max}(u_0))\times [0,1]$ and 
$w_0=\theta_0(u_0)$. We know that $w$ is radial, so for all 
$(s,y)\in (0,\frac{T_{max}(u_0)}{N^2})\times \overline{B} $, $w(s,y)=\tilde{w}(s,
\|y\|)$ with 
\begin{equation}
\label{fomule_w}
\tilde{w}\in C^{1,2}( (0,\frac{T_{max}(u_0)}{N^2})\times [0,1])
\end{equation}
We have shown that $u(t,x)=x \, \tilde{w}(\frac{t}{N^2},x^{\frac{1}{N}})$ so
that 
\begin{equation}
\label{formule_u_x}
u_x(t,x)=\tilde{w}(\frac{t}{N^2},x^{\frac{1}{N}})+\frac{x^{\frac{1}{N}}}
{N}\tilde{w}_r(\frac{t}{N^2},x^{\frac{1}{N}}) 
\end{equation}
This formula already allowed us to prove that $u(t)\in C^1([0,1])$ with 
$u_x(t,0)=\tilde{w}(\frac{t}{N^2},0)$. Since $w(t)$ is radial, then 
$\tilde{w}_r(\frac{t}{N^2},0)=0$ so we get that
$$ |u_x(t,x)-u_x(t,0)|\leq K\, x^{\frac{2}{N}} $$
with $K= (\frac{1}{2}+\frac{1}{N})\|\tilde{w}(\frac{t}{N^2})_{rr}\|_{\infty,[0,1]} $. 
Hence, $u(t)\in Y_m^{1,\frac{2}{N}}$.\\

\subsection{Shape of the derivative of classical solutions of problem $(PDE_m)$ }
\label{forme}

We will prove Proposition \ref{propo_forme}.\\

\textit{Proof :} i) We set $h(t,x)=\tilde{w}(\frac{t}{N^2},x)+\frac{x}
{N}\tilde{w}_r(\frac{t}{N^2},x)$.\\
The result comes from formula (\ref{formule_u_x}) and because of 
(\ref{fomule_w}).\\
ii) Since $\tilde{w}(\frac{t}{N^2},0)>0$ for $t\in [t_0,T]$, then by 
compactness, there exists $\delta>0$ such that 
$\tilde{w}(\frac{t}{N^2},x^{\frac{1}{N}})+\frac{x^{\frac{1}{N}}}{N}\tilde{w}_r(\frac{t}
{N^2},x^{\frac{1}{N}}) >0$ on $[t_0,T]\times [0,\delta]$. \\
Since $x\mapsto x^q$ is smooth on $(0,\infty)$ and $w$ satisfies
$$w_t=\Delta w+N^2 \, w(w+\frac{y.\nabla w}{N})^q $$
then by classical regularity result, $w\in C^{1,\infty}([t_0,T]\times \overline{B}(0,
\delta))$. This gives the regularity of $h$.\\
iii) Clear since $\tilde{w}$ has odd order derivatives vanishing at $x=0$.\\

{\bf Acknowledgements :} the author would like to thank Philippe Souplet 
for all stimulating discussions and comments about this paper.

\end{document}